\documentclass[oneside,english,11pt]{amsart}
\usepackage{}
\usepackage{txfonts}
\usepackage{amsfonts}
\usepackage{lmodern}

\usepackage[T1]{fontenc}
\usepackage[latin9]{inputenc}
\usepackage{geometry}
\usepackage{amsthm}
\usepackage{amssymb}
\usepackage{amssymb,amsmath,amsthm,tikz,tikz-cd}
\usetikzlibrary{matrix,arrows,decorations.pathmorphing}
\usepackage{appendix}
\usepackage{tikz-cd}
\usepackage{comment}

\usepackage[bookmarks=true, bookmarksopen=true,%
bookmarksdepth=3,bookmarksopenlevel=2,%
colorlinks=true,%
linkcolor=blue,%
citecolor=blue,%
filecolor=blue,%
menucolor=blue,%
urlcolor=blue]{hyperref}

\title{Augmentations and sheaves for links}
\author{Honghao Gao}

\address{Department of Mathematics \\ Michigan State University \\ C212 Wells Hall, 619 Red Cedar Road, East Lansing, MI, 48824.}
\email{gaohongh@msu.edu}

\makeatletter

\numberwithin{equation}{section}
\numberwithin{figure}{section}

\theoremstyle{plain}
\newtheorem{thm}{Theorem}[section]
\newtheorem{lem}[thm]{Lemma}

\newtheorem{prop}[thm]{Proposition}
\newtheorem{thm-defn}[thm]{Theorem-Definition}

\newtheorem{prop-def}[thm]{Proposition-Definition}

\theoremstyle{definition}
\newtheorem{defn}[thm]{Definition}

\newtheorem{notation}[thm]{Notation}

\theoremstyle{remark}
\newtheorem{rmk}[thm]{Remark}

\newcommand{\bbC}{{\mathbb{C}}}

\newcommand{\bbR}{{\mathbb{R}}}

\newcommand{\bbZ}{{\mathbb{Z}}}

\newcommand{\cA}{{\mathcal{A}}}

\newcommand{\cE}{{\mathcal{E}}}
\newcommand{\cF}{{\mathcal{F}}}
\newcommand{\cG}{{\mathcal{G}}}
\newcommand{\cH}{{\mathcal{H}}}
\newcommand{\cI}{{\mathcal{I}}}

\newcommand{\cL}{{\mathcal{L}}}
\newcommand{\cM}{{\mathcal{M}}}

\newcommand{\del}{{\partial}}

\newcommand{\la}{{\langle}}
\newcommand{\ra}{{\rangle}}
\newcommand{\ep}{{\epsilon}}

\newcommand{\idV}{{\mathrm{id}_V}}
\newcommand{\id}{{\mathrm{id}}}

\tikzset{node distance=1.5cm, auto}

\usepackage{babel}

\begin{document}
\maketitle

\begin{abstract}
We study the relation between augmentations and sheaves in the context of framed oriented links. In this set up, we find slightly more sheaves than augmentations. After removing the sporadic sheaves, we construct a bijective correspondence between augmentations of the framed cord algebra and simple sheaves that are micro-supported along the conormal bundle of the link. 
\end{abstract}

\section{Introduction}

The augmentation sheaf correspondence is an important subject in the rigidity of contact topology. The correspondence relates two invariants of Legendrian submanifolds that are development independently. One invariant is the Legendrian contact differential graded algebra (dga) and its augmentations, a Floer theoretic invariant introduced in \cite{Che}, and then generalized to broader settings \cite{EES1,EES2,EES3}. The other invariant comes from microlocal sheaf theory, whose foundation was found in \cite{KS}, and enters contact topology more recently \cite{GKS, STZ}. The first correspondence between the two was established in \cite{NRSSZ}, between a positive augmentation category and a category of simple sheaves for Legendrian links in the standard tight contact three-fold. The correspondence has been generalized to other setups along several directions. 

This paper focuses on the set up of knot contact homology, where we consider Legendrian conormal tori, that arises from oriented framed links, in the cosphere bundle of $\bbR^3$ or $S^3$. Knot contact homology was introduced in \cite{Ng1, Ng2, Ng3} with a knot theoretic formulation, and proven to be equivalent to the Legendrian dga in \cite{EENS}. As far as the augmentations concern, we can work with the framed cord algebra \cite{Ng3}, because it is the degree zero homology of the dga. On the sheaf side, we consider objects in the dg derived category of sheaves (in $\bbR^3$ or $S^3$) that are simple along the conormal bundle of the link. In this paper, we build explicit connections between them.

Let $(L,L')$ be a framed link in $\bbR^3$ or $S^3$. Let $k$ be a field. We consider the framed cord algebra of the framed link $\mathrm{Cord}(L,L')$ and its augmentations, which are algebra morphisms $\ep: \mathrm{Cord}(L,L')\rightarrow k$. An algebraic torus, whose dimension equals the number of components of $L$, acts on the set of augmentations, and we define $\cA ug$ to be the moduli set of augmentations up to this action, see Definition \ref{augmoduli}. On the sheaf side, we consider simple sheaves of microlocal rank $1$ which are micro-supported along the conormal bundle of $L$. Let $\widetilde{\cM}$ be the moduli set of such sheaves up to locally systems, modulo the the natural isomorphism in the quotient category, see \eqref{modshtil}. We observe that objects in the sheaf moduli can be expressed as an extension between a sheaf micro-supported on a sublink, and a sheaf supported on the complement link (Lemma \ref{lemShclassify}). We consider the subset $\cM\subset \widetilde{\cM}$ of sheaves arising from trivial extension classes, see Definition \ref{shmoduli}. Our main results is stated as follows.

\begin{thm}\label{Mainthm} There is a bijective correspondence
$$\cA ug \cong \cM.$$
\end{thm}

\begin{rmk}\label{Extrash}
It is a matter of perspectives whether this result should be interpreted as an ``correspondence'' or not. In the currently existing results on the augmentation-sheaf correspondence, including the special case when $L$ is a knot with the set up of this paper, the moduli space of sheaves is taken to be $\widetilde{\cM}$, instead of $\cM$ in Theorem \ref{Mainthm}. A detailed analysis of the differences between $\cM$ and $\widetilde{\cM}$ is included in Section \ref{Sec:shmod}. For many links, the inclusion $\cM\subset \widetilde{\cM}$ is strict. From this perspective, we say there are more sheaves than augmentations.

On the other hand, once we remove the sporadic objects (those in $\widetilde{\cM} \setminus\cM$) from $\widetilde{\cM}$, there is a bijective correspondence as stated in the theorem. The size of $\widetilde{\cM} \setminus\cM$ is relatively smaller than $\cM$, which can be seen from the augmentation side. To be more concrete, each object in $\widetilde{\cM} \setminus\cM$ arises from a nontrivial extension class in a certain way (see Lemma \ref{lemShclassify} for details), and it has a ``sibling'' in $\cM$ where we take the trivial extension class instead. For generic augmentations (see Definition \ref{defindexsets}), which form an Zariski open subset in $\cA ug$, their counterparts in $\cM$ do not have siblings in $\widetilde{\cM} \setminus\cM$. Since we have a complete characterization of sheaves in $\widetilde{\cM} \setminus\cM$, it would also make sense to single them out, so as to focus on the majority cases and declare the correspondence.



\end{rmk}

The correspondence in Theorem \ref{Mainthm} has been formulated for knots \cite{Gao2}. Below is a technical summary of the new features in the construction for links.

\begin{enumerate}
\item The augmentation moduli space for knots is the augmentation variety, while that for links is the quotient by an algebraic torus. See Remark \ref{augmodkvl}.

\item In the case of knots, objects in the desired sheaf category, if it is concentrated in cohomological degree $0$, can be uniquely represented by its ``stabilization'' (Definition \ref{stabilizationdefn}). The same property does not hold for links, because there could exist a nontrivial extension between sheaves that are micro-supposed along different sublinks.  See Section \ref{Sec:stsh} and Remark \ref{rsnsl}.

\item When constructing the map from sheaves to augmentations, one needs to choose a set of ``local trivialization'' functions, see Section \ref{Sec:fsta}. Different choices define equivalent augmentations in the augmentation moduli. For pure cords, the induced augmented values does not depends on the choices, and it agrees with the formula for knots. See Proposition \ref{purecordformula} and Remark \ref{fpsak}.

\item When constructing the map from augmentations to sheaves, one needs to choose a braid representative of the link, but the induced sheaf does not dependent on the braid. See Section \ref{Sec:fats}.
\end{enumerate}

Sheaves in our setting can be reduced to a representation of the link group (i.e. the fundamental group of the link complement) with extra data. Prior to the first statement of the correspondence in \cite{NRSSZ}, some of these knot group representations have been related to augmentations \cite{Ng4, Cor2}. Their results have been completed to the correspondence for knots \cite{Gao2}.

Other cases for Legendrian surfaces include: in \cite[Appendix]{CM}, Sackel proved the correspondence for Legendrian surfaces arising from planar cubic graphs, based on the dga constructed in \cite{CM}; in \cite{RS4}, Rutherford-Sullivan constructed a one-way map from augmentations to sheaves for Legendrian surfaces (with mild front singularities) in the one-jet space, based on the cellular dga constructed in \cite{RS1, RS2, RS3}. Note that in all these cases, the results are stated on the object level, because the positive augmentation category has not been rigorously defined for Legendrian surfaces (Remark \ref{nocatsta}).

We explain the intuition behind the main theorem from the perspective of contact and symplectic topology. Let $(L,L')\subset \bbR^3$ be the oriented framed link. Let $\Lambda_L\subset T^\infty \bbR^3 = \del_\infty(T^*\bbR^3)$ be the Legendrian conormal tori. Via the contact transformation $T^\infty \bbR^3 \cong \bbR^3\times S^2 = S^2 \times \bbR^3 \cong T^*S^2 \times \bbR \cong J^1(S^2)$, the Legendrian $\Lambda_L$ can be placed in the one-jet space $J^1(S^2)$. Following the recipe of \cite{Eli, EGH}, the Legendrian dga can be defined via counting holomorphic curves (in the form of gradient flow trees \cite{Ekh}) in the symplectization. The symplectization $J^1S^2\times \bbR$ can be embedded in a cotangent bundle via $J^1S^2\times \bbR \cong T^{*,+}(S^2\times \bbR) \subset T^*(S^2\times \bbR)$. The geometric counterpart for augmentations consists of Lagrangian fillings and their triangulated envelopes in the Fukaya category $Fuk(T^*(S^2\times \bbR))$, which correspond to sheaves on $S^2\times \bbR$ via microlocalization \cite{NZ, Nad}. This theoretical relation is depicted by the dashed arrow in the diagram below.

\begin{figure}[h!]
\begin{center}
\begin{tikzpicture}[
	block/.style = {draw, align=center,minimum height=3em, inner sep=5},
	arrow/.style = {->, semithick, rounded corners = 5pt}]
\node (R1) at (-5, 0) {{Floer theory}};
\node (R2) at (-5, -2.4) {{Sheaf theory}};
\node (C1) at (0, 1.2) {{$\bbR^3$}};
\node (C2) at (7, 1.2) {{$J^0S^2 = S^2\times \bbR$}};
\node[block] (Augfca) at (0, 0) {{Augmentations of} \\ {the framed cord algebra}};
\node[block] (Augdga) at (7, 0) {{Augmentations of} \\ {the Legendrain dga}};
\node[block] (Sh1) at (0, -2.5) {{Simple sheaves on $\bbR^3$,} \\ {micro-supported along $\Lambda_L$}};
\node[block] (Sh2) at (7, -2.5) {{Simple sheaves on $J^0S^2$} \\ {micro-supported along $\Lambda_L$}};	
\draw[<->] (Augfca.east) --node[above]{\cite{Ng3}}  (Augdga.west);
\node at (-1.3, -1.25) {this paper};
\draw[<->] (Sh1.north) --  (Augfca.south) ;
\draw[<->] (Sh2.west) -- node[above]{\cite{Gao1}} (Sh1.east);
\draw[<->,dashed] (Augdga.south) -- (Sh2.north);
\end{tikzpicture}
\end{center}
	\label{blueprint}
\end{figure}

Our strategy for the proof is to transform the correspondence from $T^*({S^2}\times \bbR)$ to $T^*\bbR^3$. Because the dga for $\Lambda_0$ is concentrated in non-negative degrees, augmentations of the dga is the same as augmentations of the degree $0$ homology, which is the framed cord algebra introduced in \cite{Ng3}. On the sheaf side, there is an equivalent of categories $Sh_{\Lambda_L}(\bbR^3) / Loc(\bbR^3)\cong Sh_{\Lambda_L}(S^2\times\bbR)/Loc(S^2\times \bbR)$ which preserves the simpleness of sheaves when transforming back and forth \cite{Gao1}. With both augmentations and sheaves transformed to $\bbR^3$, this paper establishes a correspondence between them.

A final remark on the diagram: the one-way construction in \cite{RS4}, if the coefficients were further enhanced to include $H_1(\Lambda_L)$, could be a candidate for the downward dashed arrow, but currently we don't know whether their map is injective, or surjective, or it commutes with the other arrows in the diagram.

\medskip
We mention some literature that are related to the context of this paper. The augmentation-sheaf correspondence has also been established in some other set-ups, such as \cite{CNS, ABS}. Augmentations can be categorified from other perspectives \cite{BC, CDRGGn, GPS, EL}. Both Floer theory and sheaf theory produce complete invariants for knots \cite{She, ENS}. The representations of the link group induced from the sheaves are related to link polynomials and character varieties, such as in \cite{Cor1, AENV}, and are potentially related to \cite{CCGLS, DG, BZ, NiZh, CS, KM, AH, GTZ, GW, HMP, MP} etc.

\medskip
This paper is organized as follows. In Section \ref{Sec:aug}, we discuss the framed cord algebra and its augmentations, and define the augmentation moduli. In Section \ref{Sec:sh}, we study properties of the sheaves, and characterize objects in the sheaf moduli. In Section \ref{Sec:correspondence}, we prove the main theorem.

\begin{notation}
Throughout the paper, we fix the following notations.
\begin{itemize}
\item
$X = S^3$ or $\bbR^3$. 
\item
$(L,L')$ is an $r$-component framed oriented link in $X$. $L = K_1\sqcup \dotsb \sqcup K_r$, and $L' = \ell_1\sqcup \dotsb \sqcup \ell_r$.
\item
For any $1\leq s\leq r$, $i_s:K_s \rightarrow X$ is the closed embedding of the component $K_s$, and $j_s: X\setminus K_s\rightarrow X$ is the open embedding of the complement.

\item
Fix a commutative ground field $k$ for representations, sheaves, and augmentations, (but not the dga).
\end{itemize}
\end{notation}

\smallskip
\noindent\textbf{Acknowledgements.} We thank Stéphane Guillermou, Efstratia Kalfagianni, Lenhard Ng, Linhui Shen and Eric Zaslow for helpful discussions and valuable comments. We thank Stéphane Guillermou for crucial help with the results in Section \ref{Sec:shmod}. This work is partially supported by ANR-15-CE40-0007 ``MICROLOCAL'' and an AMS-Simons travel grant.

\section{Augmentations}\label{Sec:aug}

\subsection{Framed cord algebra and its augmentations}

The cord algebra first appeared in \cite{Ng1, Ng2}. The framed version was introduced in \cite{Ng3}, which models the degree zero knot contact homology \cite{EENS}. We shall use a mild generalization of the version in {\cite[Definition 2.5]{CELN}}. For a review of the comparison between different versions and properties of the framed cord algebra, see \cite{Gao3}.

Let $(L,L')$ be an $r$-component framed oriented link. We decorate each $\ell_s$ with a marked point $\ast_s\in \ell_s$. Let $\ast := \{\ast_1,\dotsb, \ast_r\}$.

For two paths $c_1$ and $c_2$, we use $c_1\cdot c_2$ to denote their concatenation.

\begin{defn}
A framed cord of $(L,L')$ is a continuous map $c: [0,1]\rightarrow X\setminus L$ such that $c(0), c(1) \in L' \setminus \ast$. Two framed cords are homotopic if they are homotopic through framed cords. We write $[c]$ for the homotopy class of the cord $c$. 

A framed cord from $\ell_s\setminus \ast_s$ to  $\ell_t\setminus \ast_t$ is often denoted by $c_{st}$, and we simply say it is a cord from $K_s$ to $K_t$.

Define a non-commutative unital ring $\cA$ as follows: as a ring, $\cA$ is freely generated by homotopy classes of framed cords and extra generators $\lambda_s^{\pm 1}, \mu_s^{\pm 1}$, $1\leq s\leq r$, modulo the relations
$$
\lambda_s \cdot \lambda_s^{-1} = \lambda_s^{-1}\cdot \lambda_s = \mu_s \cdot \mu_s^{-1} = \mu_s^{-1}\cdot \mu_s = 1,\quad \textrm{ for } 1\leq s\leq r,
$$
and
$$
\begin{cases}
\lambda_s \cdot \mu_t = \mu_t\cdot \lambda_s \\
\lambda_s \cdot \lambda_t = \lambda_t \cdot \lambda_s \quad \textrm{ for } 1\leq s,t\leq r.\\
\mu_s \cdot \mu_t = \mu_t \cdot \mu_s
\end{cases}
$$
Thus $\cA$ is generated as a $\bbZ$-module by non-commutative words in homotopy classes of cords and powers of $\lambda_s$ and $\mu_t$. The powers of the $\lambda_s$ and $\mu_s$ commute with each other, but do not commute with any cords.

The framed cord algebra is the quotient ring 
$$\textrm{Cord}(L, L') = \cA/\cI,$$
where $\cI$ is the two-sided ideal of $\cA$ generated by the following relations:
\begin{itemize}
\item
(normalization) \quad $[e_s] = 1-\mu_s$, where $e_s$ is a constant cord on $K_s$,
\item
(meridian) \quad\quad\;\;\; $[m_s \cdot c_{st}] = \mu_{\{s\}}[c_{st}],\quad  [c_{st} \cdot m_t ] = [c_{st}]\mu_{\{t\}}$,
\item
(longitude) \qquad \quad $[\ell_{\{s\}}\cdot c_{st}] = \lambda_{\{s\}} [c_{st}], \quad [c_{st}\cdot \ell_{\{t\}}] = [c_{st}]\lambda_{\{t\}}$,
\item
(skein relations) \quad $[c_{sk} \cdot c_{kt}] = [c_{sk} \cdot m_k\cdot c_{kt}] + [c_{sk}][c_{kt}]$, where $c_{sk}$ and $c_{kt}$ are composable cords, and $m_k$ is a meridian based at the composing point.
\end{itemize}
\end{defn}

Up to $\bbZ$-algebra isomorphisms, $\mathrm{Cord}(L,L')$ does not depend on the decoration or the framing. Hence we can assume $L'$ is the Seifert framing, and simply denote the framed cord algebra as $\mathrm{Cord}(L)$.

Some constructions in this paper replies on a braid representative of $L$. We discuss the framed cord algebra in this setting. Suppose $L$ is the closures of an $n$-stranded braid $B \in Br_n$. Let $D$ be a disk cutting the braid transversely. Let $x_i = (i,-1)$ and $y_i=(i,0)$, and we assume $L\cup D = \{y_1, \dotsb, y_n\}$, and $L'\cup D = \{x_1,\dotsb, x_n\}$. We also assume that marked points $\ast$ are not contained in $D$.

Let $\gamma_{ij}$, $1\leq i,j\leq n$ be the linear path connecting $x_i$ and $x_j$, which is called a \textit{standard cord}. Standard cords generate the framed cord algebra \cite[Lemma 2.12]{Gao3}. Let $m_i^{(k)}$ be the meridian loop $m_i$ based at $x_k$. Relations with respect to standard cords are:
\begin{itemize}
\item
$[\gamma_{ii}] = 1 -\mu_{\{i\}},$
\item
$[m_i^{(i)}\cdot \gamma_{ij}] = \mu_{\{i\}}[\gamma_{ij}],\quad  [\gamma_{ij} \cdot m_j^{(j)}] = [\gamma_{ij}]\mu_{\{j\}},$
\item
$[\ell^{(i)}_{\{i\}}\cdot \gamma_{ij}] = \lambda_{\{i\}} [\gamma_{ij}], \quad [\gamma_{ij}\cdot \ell^{(j)}_{\{j\}}] = [\gamma_{ij}]\lambda_{\{j\}},$
\item
$
[\gamma_{it} \cdot \gamma_{tj}] = [\gamma_{it} \cdot m_t^{(t)}\cdot \gamma_{tj}] + [\gamma_{it}][\gamma_{tj}].
$
\end{itemize}

\begin{figure}[h]
	\centering
	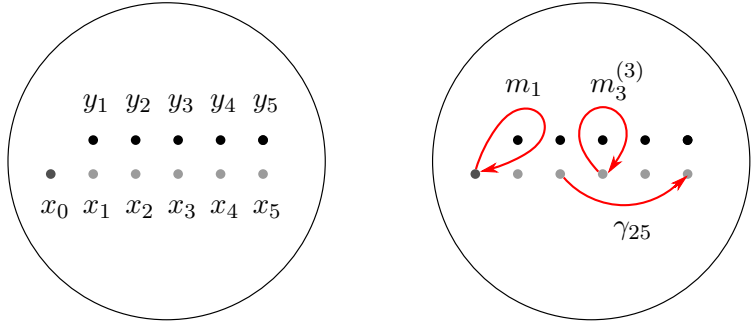
	\caption[]{Examples of meridians and framed cords on the transverse disk.}
	\label{Fig:ConfigurationDisk}
\end{figure}

Let $\{ - \}: \{1,\dotsb, n\}\rightarrow \{1,\dotsb, r\}$ be the \textit{component function}, i.e. strand $i$ belongs to the component $K_{\{i\}}$. After a suitable conjugation of $B$, we can assume that the component function is non-decreasing.

\begin{defn}
An \textit{augmentation} $\epsilon$ of $\textrm{Cord}(L)$ is a unit preserving algebra morphism
$$\epsilon: \textrm{Cord}(L)\rightarrow k,$$
where $k$ is any commutative field.
\end{defn}

For a framed cord $[c]$, we abbreviate $\epsilon([c])$ as $\epsilon(c)$. For a standard cord $\gamma_{ij}$, we simply write $\ep_{ij}$ for $\ep(\gamma_{ij})$. Suppose $L = \langle B \rangle$, where $B\in Br_n$, we define an $n\times n$ matrix $R$ by setting $R_{ij} = \ep_{ij}$. Namely,
$$
R = 
\begin{pmatrix}
\epsilon_{11} & \dotsb & \epsilon_{1n} \\
\vdots & \ddots &\vdots \\
\epsilon_{n1} & \dotsb & \epsilon_{nn}
\end{pmatrix}.
$$
Let $R_j$ be the $j$-th column vector and $R^i$ be the $i$-th row vector.

\begin{lem}
The following are equivalent:
\begin{enumerate}
\item $R_j =0$;
\item $R_k =0$ for all $j$ with $\{k\} = \{j\}$;
\item $\ep(c_{st}) = 0$ for all framed cords $c_{st}$ with $t = \{j\}$;
\end{enumerate}
\end{lem}

\begin{proof}
Clearly $(3)\Rightarrow (2)\Rightarrow (1)$.

$(1)\Rightarrow (3)$. We can homotope the framed cord $c_{st}$ such that it takes the form $\gamma_{ij} \cdot c$ where $j$ is given in (1), $\{i\}= s$ and $c$ is a loop based at $x_j$. Note $c$ can be generated by $\{m_k\}_{1\leq k\leq n}$. We proceed with an induction on the word length of $c$ with respect to the generating meridians. The base case, where $c$ is the constant loop, is trivial. Suppose $c = c' \cdot m_k^{\pm 1}$.

If $c = c' \cdot m_k$, then 
$$
\epsilon(\gamma_{ij} \cdot c' \cdot m_k)
	= \epsilon(\gamma_{ij} \cdot c' \cdot \gamma_{jk} \cdot m_k^{(k)}\cdot \gamma_{kj}) 
	= \epsilon(\gamma_{ij} \cdot c') - \epsilon(\gamma_{ij} \cdot c' \cdot \gamma_{jk}) \epsilon(\gamma_{kj})=0.
$$
The case $c = c' \cdot m_k^{-1}$ is similar:
$$\epsilon(\gamma_{ij} \cdot c' \cdot m_k^{-1}) = \epsilon(\gamma_{ij} \cdot c') + \epsilon(\gamma_{ik} \cdot \tilde{c}' \cdot m_k^{-1}) \epsilon(\gamma_{kj}) =0.$$
\end{proof}

There are similar statements with respect to the row vectors of $R$.
\begin{lem}
The following are equivalent:
\begin{enumerate}
\item $R^{i} =0$;
\item $R^{k} =0$ for all $j$ with $\{k\} = \{i\}$;
\item $\ep(c_{st}) = 0$ for all framed cords $c_{st}$ with $s = \{i\}$;
\end{enumerate}
\end{lem}

\begin{defn}\label{defindexsets} For $I = \{1,\dotsb, n\}$, 
define the partition $I = I'\cup I''$ such that $R^i \neq 0$ for $i\in I'$ and $R^i = 0$ for $i \in I''$, and the partition $I = J' \cup J''$ such that $R_j \neq 0$ for $j\in J'$ and $R_j = 0$ for $j\in J''$.

The $(i,j)$-entry of $R$ is generic if both $i, j \in I'\cap J'$.

An augmentation $\epsilon$ is \textit{generic} if $I' = J' = I$, or equivalently $I'' = J'' = \emptyset$.
\end{defn}

If particular, if $\ep(\mu_i)\neq 1$ for all $i$, then $\epsilon$ is generic because diagonal entries of $R$ are nonzero.

\subsection{Augmentation moduli}
Let $\cA ug_{naive}$ be the set of augmentations of $\mathrm{Cord}(L)$. It admits the structure of an affine variety by considering the evaluation on standard cords. Suppose $L =\la B\ra$ as the closure of an $n$-strand braid. The standard cords $\{\gamma_{ij}\}_{1\leq i,j\leq n}$ and $\{\lambda_i^{\pm 1}, \mu_i^{\pm 1}\}_{1\leq i\leq r}$ generate $\textrm{Cord}(L)$, for example see \cite[Lemma 2.12]{Gao3}. Hence,
$$\cA ug_{naive} = \left\{\big(\epsilon(\lambda_i), \epsilon(\mu_i), \epsilon(\gamma_{ij})\big) \,|\, \epsilon \textrm{ is an augmentation}\right\} \subset (k^{*})^{2r}\times (k)^{n^2}.$$
Here $\epsilon(\lambda_i), \epsilon(\mu_i)\in k^*$ because $\lambda_i,\mu_i$ are invertible in $\mathrm{Cord}(L)$. In particular, $\cA ug_{naive} \subset (k^{*})^{2r}\times (k)^{n^2} $ is the subvariety cut off by the following relations:

\begin{itemize}
\item
$\epsilon(\gamma_{ii}) = 1 - \epsilon(\mu_{\{i\}}),$

\item
$\epsilon(m_i^{(i)}\cdot \gamma_{ij}) = \epsilon(\mu_{\{i\}})\epsilon(\gamma_{ij}),\quad \epsilon(\gamma_{ij} \cdot m_j^{(j)})= \epsilon(\gamma_{ij})\epsilon(\mu_{\{j\}}),$
\item
$\epsilon(\ell^{(i)}_{\{i\}}\cdot \gamma_{ij}) = \epsilon(\lambda_{\{i\}}) \epsilon(\gamma_{ij}), \quad \epsilon(\gamma_{ij}\cdot \ell^{(j)}_{\{j\}}) = \epsilon(\gamma_{ij})\epsilon(\lambda_{\{j\}}),$

\item
$
\epsilon (\gamma_{it} \cdot \gamma_{tj}) = \epsilon(\gamma_{it} \cdot m_t^{(t)}\cdot \gamma_{tj}) + \epsilon(\gamma_{it})\epsilon(\gamma_{tj}).
$
\end{itemize}

\bigskip
However, to relate augmentations to sheaf, we need to consider a further quotient $\cA ug$.

\begin{defn}\label{defndilation}
A dilation parameter for $L$ is an $r$-tuple $d = (d_1,\dotsb, d_r)\in (k^*)^r$.
\end{defn}

\begin{lem}\label{dilationaction}
Let $(L,L')$ be an $r$-component framed oriented link. Let $d = (d_1,\dotsb, d_r)\in (k^*)^r$ be a dilation parameter.

For any augmentation $\epsilon: \mathrm{Cord}(L)\rightarrow k$, we define a map $\epsilon':  \mathrm{Cord}(L)\rightarrow k$ by 
$$\epsilon'(\lambda_s) = \epsilon(\lambda_s),\quad \epsilon'(\mu_s) =\epsilon(\mu_s), \quad \epsilon'(c_{st}) = \frac{d_s}{d_t}\,\epsilon(c_{st}),$$
where $c_{st}$ is a framed cord from $K_s$ to $K_t$. Then $\epsilon'$ is also an augmentation.
\end{lem}

\begin{proof}
We check the relations in the framed cord algebra. For the normalization, suppose $e_s$ is the constant cord on $K_s$, then $\epsilon'(e_{s}) = \epsilon(e_{s}) = 1-\epsilon(\mu_s) = 1-\epsilon'(\mu_s).$ Meridian and longitude relations follow in a similar way. If $c_{sk}$, $c_{kt}$ are two composable paths, the skein relation yields $[c_{sk} \cdot c_{kt}] = [c_{sk} \cdot m_k^{(k)}\cdot c_{kt}] + [c_{sk}][c_{kt}].$ Then
\begin{align*}
\epsilon'(c_{sk} \cdot c_{kt}) &= \frac{d_s}{d_t}\,\epsilon(c_{sk} \cdot c_{kt}) \\
	&=\frac{d_s}{d_t}\, \big(\epsilon(c_{sk} \cdot m_k^{(k)}\cdot c_{kt}) + \epsilon(c_{sk})\epsilon(c_{kt})\big) \\
	&= \frac{d_s}{d_t}\, \big(\frac{d_t}{d_s}\, \epsilon'(c_{sk} \cdot m_k^{(k)}\cdot c_{kt}) + \frac{d_k}{d_s}\,\epsilon'(c_{sk})\frac{d_t}{d_k}\,\epsilon'(c_{kt})\big) \\
	&= \epsilon'(c_{sk} \cdot m_k^{(k)}\cdot c_{kt}) + \epsilon'(c_{sk})\epsilon'(c_{kt}).
\end{align*}
\end{proof}

Let $T:= (k^*)^r$, then it follows from the lemma that $T$ acts on $\cA ug_{naive}$ by dilations. 

\begin{defn}\label{augmoduli}
We define the moduli set of augmentations as
$$\cA ug := \cA ug_{naive}/ T.$$
\end{defn} 

\begin{rmk}\label{augmodkvl}
The dilation action is not free, because $(d_1,d_2, \dotsb, d_n)$ and $(1,d_2/d_1, \dotsb, d_n/d_1)$ send an augmentations $\epsilon$ to the same augmentation $\epsilon'$. We say $d$ is reduced if $d_1 = 1$. Then $$\cA ug = \cA ug_{naive}/ \{\textrm{reduced dilations}\}.$$

\noindent Reduced dilations are trivial for knots, i.e. $\cA ug = \cA ug_{naive}$. Hence, dilations are not considered in the set up of \cite{Gao2}.
\end{rmk}

\begin{rmk}\label{nocatsta}
The notion of the dilation action is borrowed from \cite[Sec. 5.3]{NRSSZ}. As far as the author knows, the positive augmentation category, $\cA ug_+(\Lambda)$, has not been rigorously defined for Legendrian surfaces. Our ad-hoc Definition \ref{defndilation} mimics the isomorphism in $\cA ug_+(\Lambda)$ for Legendrian links. In principle, the isomorphism consists of the contribution from degree $-1$ Reeb chords and weights of multi-components. There are no such Reeb chords in knot contact homology, and the remaining factor is the action of dilations as stated.
\end{rmk}

\section{Sheaves}\label{Sec:sh}

We discuss microlocal sheaves in the context of conormal of links. In Section \ref{Sec:shintro}, we review the basics of microlocal sheaf theory. In Sections \ref{Sec:shred} and \ref{Sec:stsh}, we study some properties that are useful in building the correspondence. In Section \ref{Sec:shmod}, we define the sheaf moduli spaces and describe representatives in those spaces.

\subsection{Invariants}\label{Sec:shintro}

We first recall some facts from microlocal sheaf theory. See \cite{KS, Gui} for more details.

Let $Y$ be a smooth manifold. Let $Mod(Y)$ be the abelian category of sheaves of $k$ modules on $Y$, and $Sh(Y)$ be its dg derived category. Any object $\cF\in Sh(Y)$ defines a closed conic involutive subset $SS(\cF)\subset T^*X$ called its micro-support. Let $T^\infty Y := \del_\infty (T^*Y)$ be the contact cosphere bundle. Fix a smooth Legendrian $\Lambda \subset T^\infty Y$, the full dg subcategory $Sh_\Lambda(X) \subset Sh(Y)$
$$Sh_\Lambda(Y) := \{\cF\in Sh(Y) \,|\, SS(\cF) \cap T^\infty Y \subset \Lambda\},$$
is a Legendrian isotopy invariant of $Y$, following the main theorem of \cite{GKS}. Locally constant sheaves form a full subcategory, denoted by $Loc(Y)$, of $Sh_\Lambda(Y)$. We often consider the dg quotient $Sh_\Lambda(Y)/Loc(Y)$, and its variations. Let $loc(Y): = Loc (Y)\cap Mod(Y)$.

A sheaf $\cF\in Sh(Y)$ is (microlocally) simple along $\Lambda$ if $\mu hom (\cF,\cF)|_{T^\infty Y} = k_\Lambda$, where $\mu hom$ is defined in \cite[Section 4.4]{KS}. Alternatively described, for any $p = (x,\xi )\in \Lambda$, choose a local function $\phi: B_x(\epsilon) \rightarrow \bbR$ such that $\phi(x) = 0$, $d\phi_x = \xi$, then the microlocal Morse cone, defined by $\mu_{p,\phi}(\cF):= R\Gamma_{\{\phi\geq 0\}}(\cF)_x$, is isomorphic to $k[d]$ of some degree $d$. Note the degree $d$, but not the simpleness, depends on the choices $\phi$. Let $Sh^s_\Lambda(Y)\subset Sh^s(Y)$ be the full subcategory of sheaves that are simple along $\Lambda$.

When $\Lambda$ is equipped with a Maslov potential $\mu$, one can consider simple sheaves whose microlocal local Morse cones are compatible with $\mu$, assuming there is a fixed choice of local functions. Denote it by $Sh_\Lambda^{s,\mu}(X)\subset Sh_\Lambda^{s}(X)$. When $\Lambda$ is connected, there are $\bbZ$ worth of Maslov potentials, and $Sh^s_\Lambda(Y)$ is decomposed into $\bbZ$ isomorphic subcategories. Hence one could write $Sh_\Lambda^{s,\mu}(X)$ without specifying $\mu$. Whereas in the case that $\Lambda$ contains multiple components, one must specify $\mu$.

We turn to the context of this paper. Let $X = \bbR^3$ or $S^3$, and $L = K_1\sqcup \dotsb \sqcup K_r$ an $r$ component smooth link. Let $N^*_LX$ be its conormal bundle and $\Lambda_L : = N^*_LX\cap T^\infty X$ be the Legendrian tori. We consider the Maslov potential that are $0$ on all link components. Note that this is the same Maslov potential that implicitly built in the knot contact homology. For any point $p\in T^*\Lambda$, we can find a local chart $\bbR^3_{x_1x_2x_3}\cong U\subset X$ such that $$L\cap U = \{x_2=x_3 =0\},$$
and $p = dx_3$ at the origin. We choose $\phi = x_3 + x_1^2 + x_2^2$. Then our Maslov potential imposes the condition $\mu_{p,\phi}(F) = k$. Let $Sh_{\Lambda_L}^{s,0}(X)$ be the full dg subcategory of such sheaves.

Let $j: X\setminus L\rightarrow X$ be the open embedding, and $i_s: K_s\rightarrow X$ be the closed embeddings.

\begin{lem}\label{decomploc}
$F\in Sh_{\Lambda}(X)$ if and only if $j^{-1}\cF \in Loc(X\setminus L)$ and $i_s^{-1}\cF \in Loc(K_s)$. $F\in Mod_{\Lambda}(X)$ if and only if $j^{-1}\cF \in loc(X\setminus L)$ and $i_s^{-1}\cF \in loc(K_s)$. 
\end{lem}
\begin{proof}
Because the micro-support is locally defined, the proof follows \cite[Lemma 3.1]{Gao2}.
\end{proof}

Let $\pi_L := \pi_1(X\setminus L)$ be the link group. Fix a framing $L' = \ell_1\sqcup \dotsb \sqcup \ell_r$ of $L$. Let $m_s$ be a fixed meridian of $K_s$.

\begin{lem}\label{decompabsheaf}
A sheaf $\cF \in Mod_{\Lambda_L}(X)$ is equivalent to the following data:
\begin{enumerate}
\item a representation $\rho: \pi_L\rightarrow GL(V)$, and
\item for each $1\leq s\leq r$, a representation $\rho_s: \bbZ \rightarrow GL(W_s)$, and
\item a linear transformation $T_s: W_s\rightarrow V$, such that (a) $\rho(\ell_s) \circ T_s = T_s \circ \rho'(K_s)$ and (b) $m_s$ acts on the image of $T_s$ as identity.
\end{enumerate}
Moreover, if $\cF \in Sh_{\Lambda_L}^{s,0}(X)\cap Mod_{\Lambda_L}(X)$, then $T_s$ is surjective and the cokernel has rank $1$.
\end{lem}

\begin{proof}
Because the micro-support is locally defined, the proof follows \cite[Lemma 3.3]{Gao2}.
\end{proof}

We write $\cF\leftrightarrow (\rho, V, \rho_s,W_s, T_s)$ to indicate these equivalent data. Similarly, $\cE \leftrightarrow (\rho,V)$ for a locally constant sheaf $\cE\in loc(Y)$ and the corresponding representation $\rho: \pi_1(Y)\rightarrow GL(V)$.

\begin{rmk}\label{fungro}
Let $\Pi_1$ be the fundamental groupoid. The representations $\rho$ (resp. $\rho_s$) in the lemma can be regarded as a representation of $\Pi_1(X\setminus L)$, (resp. $\Pi_1(K_s)$).

More specifically, a path $c: [0,1]\rightarrow X\setminus L$ induces an isomorphism via pull back
$$c^*: \cF_{c(1)}\xrightarrow{\sim} \cF_{c(0)}.$$
Because both sides are isomorphic to $V$, we can trivialize the map to be
$$A_c: V\rightarrow V.$$
Suppose $c_1$ and $c_2$ are composable, then $A_{c_1\cdot c_2} = A_{c_1}\circ A_{c_2}$ because $(c_1\cdot c_2)^* = c_1^* \circ c_2^*$.
\end{rmk}

\begin{rmk}
Suppose $L$ is represented by an $n$-strand braid $B$. We can take a disk $D$ that is transverse to the braid. As a variation of Lemma \eqref{decompabsheaf}, a sheaf $\cF \in Sh_{\Lambda_L}^{s,0}(X)\cap Mod_{\Lambda_L}(X)$ is equivalent to the data of $(V,\rho, W_i, \rho_i, T_i)$, $1\leq i\leq n$, with constraints on the subspaces $W_i$ as follows.

Let $\tau: Br_n\rightarrow S_n, B \mapsto \tau_B$ be the natural map from the braid group to the permutation group. For $1 \leq i \leq n$, let $\ell_{i,\tau_B(i)}$ be the segment of $\ell_{\{i\}}$ connecting $x_{i}$ and $x_{\tau_B(i)}$. Similarly we define $K_{i,\tau_B(i)}$. Each $\rho_i: W_i \stackrel{\sim}{\rightarrow} W_{\tau_B(i)}$ is determined by $K_{i,\tau_B(i)}$. Then the automorphism $A_{\ell_{i,\tau_B(i)}}: V\rightarrow V$ satisfies $A_{\ell_{i,\tau_B(i)}}|_{W_{\tau_B(i)}} = \rho_i: W_{\tau_B(i)} \stackrel{\sim}{\rightarrow} W_i$. 
\end{rmk}

\subsection{Reduced sheaf}\label{Sec:shred}

\begin{defn}\label{redsheaf}
A sheaf $\cF\in Sh^{s,0}_{\Lambda_L}(X)\cap Mod(X)$ is \textit{reduced} if satisfies the following constraints.
\begin{enumerate}
\item
There does not exist an exact sequence $0\rightarrow \cL_X\rightarrow \cF \rightarrow \cF' \rightarrow 0$, with $\cL_X\in loc(X)$.

\item
There does not exist an exact sequence $0\rightarrow \cF' \rightarrow \cF \rightarrow \cL_X \rightarrow 0$, with $\cL_X\in loc(X)$.

\item
$\cF$ does not have a direct summand $\cF'$, such that $0\rightarrow \cF' \rightarrow \cL_X \rightarrow i'_*k_{L'}\rightarrow 0,$ where $L'\subset L$ is a sublink and $i': L' \rightarrow X$ is the closed embedding.
\end{enumerate}
\end{defn}

\begin{defn}
A sheaf $\cF\in Sh^{s,0}_{\Lambda_L}(X)\cap Mod(X)$ is a \textit{strongly reduced} if 
$$R^0\Gamma(\cF) =0, \quad R^0\Gamma(D'\cF) =0.$$
\end{defn}

\begin{rmk}
Recall from \cite{KS}, The dual sheaf $D'\cF$ is defined to be 
$$D'\cF := R\cH om(\cF, k_X).$$
The operator $D': Sh(X)^{op}\rightarrow Sh(X)$ is an involution, i.e. $D'(D'\cF) = \cF$. The operator is related to the Vertier dual $D: Sh(X)^{op}\rightarrow Sh(X)$. Suppose $a: X\rightarrow \textrm{pt}$ is the projection to a point, then
$$D'\cF = R\cH om(\cF, a^{-1}k_{\textrm{pt}}), \quad D\cF = R\cH om(\cF, a^{!}k_{\textrm{pt}}).$$
Since $X$ is oriented in our set up, the two dual sheaves are related by a cohomological degree shift $D'\cF = D\cF[\dim X] = D\cF[3]$.

In general, even if $\cF$ is concentrated at degree zero, $D'\cF$ may or may not be concentrated at degree zero. But if $\cL_X\in loc (X)$, then $D'\cL_X \in loc(X)$.
\end{rmk}

\begin{prop}\label{strongredeq}
An object $\cF\in Sh^{s,0}_{\Lambda_L}(X)$ is strongly reduced if and only if it satisfies all of the following constraints.
\begin{enumerate}
\item
There does not exist an exact sequence $0\rightarrow \cL_X\rightarrow \cF \rightarrow \cF' \rightarrow 0$, with $\cL_X\in loc(X)$.

\item
There does not exist an exact sequence $0\rightarrow \cF' \rightarrow \cF \rightarrow \cL_X \rightarrow 0$, with $\cL_X\in loc(X)$.

\item
$\cF$ does not have a quotient $\cF'$, such that $0\rightarrow \cF' \rightarrow \cL_X \rightarrow i'_*k_{L'}\rightarrow 0,$ for a sublink $L'\subset L$, where $i': L' \rightarrow X$ is the closed embedding, i.e. there does not exist 
$$0\rightarrow \cF'' \rightarrow \cF \rightarrow \cL_X \rightarrow i'_*k_{L'}\rightarrow 0.$$
\end{enumerate}
\end{prop}

\begin{proof}
We prove both directions by contrapositive.

 $(\Rightarrow)$ If there exists a short exact sequence $0\rightarrow \cL_X \rightarrow \cF \rightarrow \cF'\rightarrow 0$, then applying $\Gamma(-)$ we get $0 \rightarrow \Gamma(\cL_X) \rightarrow \Gamma(\cF) \rightarrow \Gamma(\cF')$. Hence, $R^0\Gamma(\cF) = \Gamma(\cF) \neq 0$.

Suppose $0\rightarrow \cF' \rightarrow \cF \rightarrow \cL_X \rightarrow 0.$ Applying $D'$, there is a distinguished triangle $D'\cL_X \rightarrow D'\cF \rightarrow D'\cF' \xrightarrow{+1}$. Since the hom functor is left exact, $D'\cF$ and $D'\cF'$ are only concentrated on non-negative degrees. Taking $R\Gamma$ yields an long exact sequence
$$0 \rightarrow R^0\Gamma(D'\cL_X) \rightarrow R^0\Gamma(D'\cF) \rightarrow R^0\Gamma(D'\cF')\rightarrow \dotsb.$$
Consequently, $R^0\Gamma(D'\cF) \neq 0$.

In the third case, we first show that $R^0\Gamma(D'i_{s*}k_{K_s}) = R^1\Gamma(D'i_{s*}k_{K_s}) =0$. Consider the short exact sequence
\begin{equation}\label{openclosedsss}
0\rightarrow j_{s!}k_{X\setminus K_s} \rightarrow k_X \rightarrow i_{s*}k_{K_s} \rightarrow 0.
\end{equation}
Since $D'k_X = k_X$, we have $R^0\Gamma(D'k_X) = k$, $R^{1}\Gamma(D'k_X) = 0$. We also have $R^0\Gamma(D'j_{s!}k_{X\setminus K_s}) = k$, because for any $i\in \bbZ$, there is
\begin{align*}
R^i \Gamma(D' j_{s!}k_{X\setminus K_s}) &= R^i \textrm{Hom} (k_X, R\cH om(j_{s!}k_{X\setminus K_s}, k_X)) \\
	&= R^i \textrm{Hom} (k_X\otimes j_{s!}k_{X\setminus K_s}, k_X) \\
	&= R^i \textrm{Hom} (k_{X\setminus K_s}, j_{s}^{-1}k_X) \\
	&= H^i (X \setminus K_s, k).
\end{align*}

Apply $R\Gamma\circ D'$ to \eqref{openclosedsss} and take the long exact sequence, we get
$$0 \rightarrow R^0\Gamma(D'i_{s*}k_{K_s}) \rightarrow R^0\Gamma(D'k_X) \rightarrow R^0\Gamma(D' j_{s!}k_{X\setminus K_s}) \rightarrow R^1\Gamma(D'i_{s*}k_{K_s}) \rightarrow R^1\Gamma(D'k_X)\rightarrow \dotsb.$$
It follows the previous calculation that $R^0\Gamma(D'i_{s*}k_{K_s}) = R^1\Gamma(D'i_{s*}k_{K_s}) =0$.

Now suppose there is a short exact sequence $0\rightarrow \cF' \rightarrow \cL_X \rightarrow i'_*k_{L'}\rightarrow 0.$ Apply $R\Gamma\circ D'$ and take the long exact sequence,
$$0 \rightarrow R^0\Gamma(D'i'_*k_{L'}) \rightarrow R^0\Gamma(D'\cL_X) \rightarrow R^0\Gamma(D'\cF') \rightarrow R^1\Gamma(D'i'_*k_{L'}) \rightarrow  \dotsb.$$
Therefore $R^0\Gamma(D'\cF') = R^0\Gamma(D'\cL_X) \neq 0$.

Apply $R\Gamma\circ D'$ to $0\rightarrow \cF'' \rightarrow \cF \rightarrow \cF'\rightarrow 0$, we get
$$0 \rightarrow R^0\Gamma(D'\cF') \rightarrow R^0\Gamma(D'\cF) \rightarrow R^0\Gamma(D'\cF'') \rightarrow  \dotsb.$$
Since $R^0\Gamma(D'\cF')\neq 0$, we have $R^0\Gamma(D'\cF)\neq 0$ as well.

\medskip
$(\Leftarrow)$ Let $a: X\rightarrow \textrm{pt}$ be the projection. 

Suppose $R^0\Gamma(\cF)\neq 0$. There is a natural embedding $0\rightarrow a^{-1}a_{*\cF} \rightarrow \cF$. The morphism is by functoriality and the embedding can be checked stalk-wise.

Suppose $R^0\Gamma(D'F)\neq 0$. Since $D'\cF$ is concentrated only on non-negative degrees, there is a natural morphism $H^0(D'\cF)\rightarrow D'\cF$. Let $\iota: a^{-1}a_*H^0(D'\cF) \rightarrow H^0(D'\cF) \rightarrow D'\cF$ be the composition, and complete it to a distinguished triangle
\begin{equation}\label{strredtri}
a^{-1}a_*H^0(D'\cF) \xrightarrow{\iota} D'\cF \rightarrow cone(\iota) \xrightarrow{+1}.
\end{equation}
Define $\cL_X := D' a^{-1}a_*H^0(D'\cF)$. The dual triangle of \eqref{strredtri} is
\begin{equation}\label{strreddualtri}
D'cone(\iota) \rightarrow \cF \xrightarrow{D'\iota} \cL_X \xrightarrow{+1}.\end{equation} 

Since $D'\cL_X = D'D'a^{-1}a_*H^0(D'\cF) = a^{-1}a_*H^0(D'\cF)\in loc(X)$, we have $\cL_X \in loc(X)$ too. Therefore $D'cone(\iota)$ is concentrated in degrees $0$ and $1$. The cohomology sequence of \eqref{strreddualtri} is:
\begin{equation}\label{minimalpropes1}
0 \rightarrow \cH^0 \rightarrow \cF \xrightarrow{D'\iota} \cL_X \rightarrow \cH^1 \rightarrow 0,
\end{equation}
where $\cH^i  := H^i(D'cone(\iota))$ for $i =0,1$. In particular, $\cH^0 \cong \ker (D'\iota)$ and $\cH^1 \cong \textrm{coker}(D'\iota)$. 

Both $\cH^0$ and $\cH^1$ are micro-supported within $\Lambda_L$, argued as follows. Let $U = X\setminus L$ and $j: U\rightarrow X$ be the open embedding. Apply the exact functor $j^{-1}$ to \eqref{minimalpropes1}, we get $0\rightarrow j^{-1}\cH^0 \rightarrow j^{-1}\cF \xrightarrow{j^{-1}D'\iota} j^{-1}\cL_X\rightarrow j^{-1}\cH^1\rightarrow 0$. Note that $j^{-1}D'\iota = D_U' j^{-1}\iota$. Since $j^{-1}$ is exact, we have
\begin{align*}
D_U'j^{-1}\cL_X &= j^{-1}a^{-1}a_*H^0(D'\cF) = a_U^{-1}a_{U*}H^0(D_U'j^{-1}\cF) = a_U^{-1}a_{U*}D_U'j^{-1}\cF.
\end{align*}
Here $a_U: U \rightarrow \textrm{pt}$, and the last equality is because $j^{-1}\cF \in loc(U)$. Therefore the morphism $D_U'j^{-1}D'\iota = D_U'D_U' j^{-1}\iota = j^{-1}\iota$ is the canonical morphism
\begin{equation}\label{minimalpropes2}
a_U^{-1}a_{U*}D_U'j^{-1}\cF \rightarrow D_U'j^{-1}\cF,
\end{equation}
and it is locally constant. Therefore $j^{-1}\cH^0 = \ker (j^{-1}D'\iota)$ is a local system. Another observation  is that \eqref{minimalpropes2} is an injection, and it further implies that $j^{-1}\cH^1 = 0$. We can argue in a similar way that $i_s^{-1}\cH^0$ is a local system on $K_s$. By Lemma \ref{decomploc}, $\cH^0$ is micro-supported within $\Lambda_L$. Finally by the triangular inequality of micro-support, we see that $\cH^1$ is also microsupported within $\Lambda_L$.

Consider \eqref{minimalpropes1}. For each component $K_s$, exactly one of $\cH^0$ or $\cH^1$ is locally constant in a neighborhood of $K_s$, and the other one is simple along $\Lambda_{K_s}$.

Suppose $\cH^1$ is locally constant near all $K_s$. Then $\cH^1 = 0$ because $j^{-1}\cH^1 =0$. We have a short exact sequence
$$0\rightarrow \cH^0 \rightarrow \cF \rightarrow \cL_X \rightarrow 0,$$
which contradicts case (2) in the assertion.

If $\cH^1$ is not locally constant near all $K_s$, let $L'\subset L$ be the sublink such that $SS(\cH^1) = \Lambda_{L'}$. Consider \eqref{minimalpropes1}, the total complex is acyclic and only $\cF,\cH^1$ are simple along points in $\Lambda_{L'}$. Since the microlocal Morse cone of $\cF$ is concentrated at degree $0$, the microlocal Morse cone of $\cH^1$ has to be at degree $1$. Combining with the fact that $j^{-1}\cH^1 =0$, we see that $\cH^1$ is isomorphic to a rank $1$ locally system supported along the link component. Since $\cL_X$ restricted to $L'$ has trivial monodromy, $\cH^1$ must also has trivial monodromy along any component of $L'$, concluding that $\cH^1 = i'_*k_{L'} $. 
We reach a contradiction to case (3) in the assertion. 
\end{proof}

\subsection{Stable sheaf}\label{Sec:stsh}

Let $\cF \in Sh^{s,0}_{\Lambda_K}(X)\cap Mod(X)$. Let $\cF \leftrightarrow (V,\rho, W_s, \rho_s, T_s)$. We choose a set of meridian generators $\{m_t\}_{t\in I}$ for $\pi_L$. For each $t \in I$, we define $V_t = \textrm{im}(\idV - \rho(m_t)) \subset V$ and define
$$V_0 := \sum_{t\in I} {V_t}\subset V.$$

Note that $V_0$ is closed under the $\pi_L$ action, because for any meridian $m_t$ and any $v \in V_0$, we have $\rho(m_t)(v) = (\rho(m_t) -\idV)(v) + \idV(v)\subset V_t + V_0 = V_0$. We denote by $(\rho_0, V_0)$ this \textit{once stabilized subrepresentation}.

\begin{lem}
$(\rho_0,V_0)$ does not depend on the choice of the generating set $\{m_t\}_{t\in I}$.
\end{lem}
\begin{proof}
Fix a meridian generating set $\{m_t\}_{t\in I}$ and $V_0$ is the associated invariant vector space. For any other meridian $m'$, $\{m_t\}_{t\in I}$ and $\{m_t\}_{t\in I}\cap\{m'\}$ define the same vector space. Suppose $m' = m_{i_1}^{\pm 1}m_{i_2}^{\pm 1}\dotsb m_{i_k}^{\pm 1}$. Let $M =\rho(m)$, keeping super and subscripts. Then for any $v\in V_0$, 
\begin{align*}
M'(v) -v &= M_{i_1}^{\pm 1}M_{i_2}^{\pm 1}\dotsb M_{i_k}^{\pm 1}(v) - v \\
	&= (M_{i_1}^{\pm 1} -\idV)M_{i_2}^{\pm 1}\dotsb M_{i_k}^{\pm 1}(v) + (M_{i_2}^{\pm 1} -\idV )\dotsb M_{i_k}^{\pm 1}(v) + (M_{i_k}^{\pm 1} - \idV)(v) \\
	&\subset V_{i_1} + V_{i_2} + \dotsb + V_{i_k} \subset V_0.
\end{align*}

If we have two generating sets of meridians, then the union is also a generating set of meridians, which defines the same $V_0$ as either one of the original generating set. 
\end{proof}

Observe that the quotient of $(\rho_0, V_0) \rightarrow (\rho,V)$ is a trivial representation. The next proposition characterizes $(\rho_0,V_0)$ as the unique, smallest sub-representation that gives a trivial quotient. (The uniqueness follows from the university property.)

\begin{prop}\label{stabsubuniv}

Any subrepresentation $(\rho',V')$ of $(\rho,V)$, such that the quotient is a trivial representation of positive dimension, contains $(\rho_0,V_0)$ as a subrepresentation.
\end{prop}
\begin{proof}
It suffices to show that $V_0\subset V'$ as a subspace. Let $\bar V = V/V'$. The group action on $\bar{V}$ is trivial, and in particular the action of any meridian is trivial. It follows from definition that $\rho(m)(v) \in V'$ for any $v\in V$, or $(\rho(m) - \textrm{id}_{V})(v) \in V'$. Since $V_0$ is spanned by the image of $(\rho(m) - \textrm{id}_{V})$ for all meridians, we conclude $V_0\subset V'$ as desired. \end{proof}

Next, we define the once stabilized subsheaf $\cF_0$. Let $(\rho_0,V_0) \leftrightarrow \cE_0\in loc(X\setminus L)$ and let $j: X\setminus L \rightarrow X$ be the open embedding. Note in the following definition $j^{-1}\cF_0 = \cE_0$, but $\cF_0$ may not be $j_*\cE_0$ when $L$ is a link.

\begin{defn}\label{stabilizationdefn}
Suppose $\cF \in Sh^{s,0}_{\Lambda_L} (X) \cap Mod(X)$ is equivalent to $(V,\rho, W_s,\rho_s, T_s)$, $1\leq s\leq r$. Define the \textit{once stabilized subsheaf} $\cF_0$ to be the sheaf equivalent to the following data $(V_0, \rho_0, W_{0s}, \rho_{0s}, T_{0s})$:
$$(\rho_0, V_0),\quad W_{0s} = W_s \cap V_0, \quad T_{0s}: W_{0s} \xrightarrow{T_s|_{W_{0s}}} V_0.$$
\end{defn}

\begin{prop}\label{F0sub}
$\cF_0$ is well-defined. $\cF_0$ is a subsheaf of $\cF$.
\end{prop}
\begin{proof}
(1) Need to check (i) $\rho_{0s}$ are well-defined and (ii) $T_{0s}$ is compactible with the corresponding peripheral subgroup action.

(i) For any $w\in W_s\cap V_0$, there is $\rho_s(K)(w) = \rho(\ell)T_s(w)\in V_0 \subset W_{0s}$. Therefore $\rho_{0s}: \bbZ\rightarrow GL(W_{0s})$ is naturally induced as a sub-representation from $\rho_s: \bbZ \rightarrow GL(W_s)$. 

(ii) $T_{0s}$ identifies $W_{0s}$ as a subspace of $V_0$. Because $(\idV - \rho(m))(w) = 0$ for any $w \in W_s$, it also holds for any $w\in W_{0s}$ --- meridian condition checked. For the longitude, the relation $T_s\circ \rho_s(K_s) = \rho(\ell_s) \circ T_s$ naturally restricts to $T_{0s}\circ \rho_{0s}(K_{s}) = \rho_0(\ell_{s}) \circ T_{0s}$ by construction.

(2) Since $W_{0s}$ is a subspace of $W_s$, $V_0$ is a subspace of $V$, together with the compatible maps, we see immediately that $\cF_0$ is a subsheaf of $\cF$.
\end{proof}

\begin{prop}
$\cF_0$ is microsupported along $\Lambda_{L_0}$ for a sublink $L_0\subset L$ ($L_0$ may or may not equal to $L$). $\cF_0$ is simple along its micro-support.
\end{prop}
\begin{proof}
By construction, $\cF_0$ restricted to the link complement or each component of the link is a local system. By Lemma \ref{decomploc},  $\cF_0$ is microsupported along a subset of $\Lambda_L$. Since $SS(\cF_0)$ is a closed subset in $T^\infty X$, there exists $L_0\subset L$ such that $SS(\cF_0)\cap T^\infty X = L_0$.

For any $p=(x,\xi)\in \Lambda_{L_0}$, there exists an open neighborhood $x\in U$ in $X$. By construction, $\cF_0|_{U}\rightarrow \cF|_{U}$ is an embedding (Proposition \ref{F0sub}) and the cone is a constant sheaf, because $W_{0s}$ is constructed by the pull back diagram $W_s \rightarrow V \leftarrow V_0$ which yields an isomorphism $V/V_0 \cong W_s/W_{0s}$. Therefore, $\mu hom(\cF_0|_U, \cF_0|_U)_p = \mu hom(\cF|_U, \cF|_U)_p  = k$, concluding that $\cF_0$ is simple along its micro-support.
\end{proof}

\begin{lem}\label{stableglobalseclem}
If $\Gamma(\cF) = 0$, then $\Gamma(\cF_0)= 0$.
\end{lem}
\begin{proof}
We prove by contrapositive. Suppose $\Gamma(\cF_0) \neq 0$. Applying $R\Gamma(-)$ to $0\rightarrow \cF_0\rightarrow \cF \rightarrow \cF/\cF_0\rightarrow 0$ and then taking the long exact sequence, we get 
$$0\rightarrow \Gamma(\cF_0) \rightarrow \Gamma(\cF) \rightarrow \Gamma(\cF/\cF_0) \rightarrow R^1\Gamma(\cF_0) \rightarrow \dotsb.$$
Since $\Gamma(\cF_0) \neq 0$, there must be $\Gamma(\cF) \neq 0$, as expected.
\end{proof}

\begin{defn}\label{stablesheafdefn}
A sheaf $\cF \in Sh^{s,0}_{\Lambda_L}(X)\cap Mod(X)$ is \textit{stable} if $\Gamma(\cF) = 0 $ and $\cF = \cF_0$.
\end{defn}

\begin{rmk}
The condition $\cF = \cF_0$ does not necessarily imply that $\Gamma(\cF) = 0$. For example, let $L = K _1\sqcup K_2$ be the two-component unlink. The link group is a free group generated by two meridians. Define $\rho:\pi_L\rightarrow GL(2,k)$ by
$$\rho(m_1) = \begin{pmatrix}
1 & 1 \\ 0 & 1
\end{pmatrix},
\quad
\rho(m_2) = \begin{pmatrix}
1 & 0 \\ 0 & \mu
\end{pmatrix},$$
where $\mu\neq 1$. Let $\cE \in loc (X\setminus L)$ be the associated local system and let $\cF = j_*\cE \in Sh_{\Lambda_L}(X)$. We see that $\cF= \cF_0$, but $\Gamma(\cF) \neq 0$ (because the section generated by $(1,0)^t$ is global).

\end{rmk}

\begin{prop}\label{sr=stab}
$\cF\in Sh^{s,0}_{\Lambda_L}(X)\cap Mod(X)$ is stable if and only if it is strongly reduced.
\end{prop}
\begin{proof}

We prove both directions by contrapositive. We assume $\Gamma(\cF) = 0$ throughout the proof.

Suppose $\cF$ is not stable, consider the short exact  sequence $0\rightarrow \cF_0 \rightarrow \cF \rightarrow \cF/\cF_0\rightarrow 0$. Denote $\bar\cF := \cF/\cF_0$, and suppose $\bar\cF \leftrightarrow (\bar\rho, \bar V)$ where $\bar V = V/V_0$. Note $ (\bar\rho, \bar V)$ is a trivial representation. If $SS(\cF_0) = SS(\cF)$, then $\bar \cF$ is a constant sheaf. By Proposition \ref{strongredeq}, $\cF$ is not strongly reduced. If $SS(\cF_0)$ is a proper subset of $SS(\cF)$, then $\bar\cF$ is micro-supported along $\Lambda_{\bar L}$ for a sublink $\bar L$. Let $j: X\setminus \bar L \rightarrow X$ be the open embedding and let $\cL_X = j_*j^{-1}\bar\cF$. Consider the natural morphism $\bar \cF \rightarrow j_*j^{-1}\bar\cF$. The morphism is an isomorphism restricted to $X\setminus \bar L$, hence the mapping cone is a sheaf only supported on $\bar L$. Because the cone has microlocal rank $1$, it is isomorphic to a direct sum of $i_{s*}\cG_s$ where $\cG_s \in loc(K_s)$ has rank $1$ and $K_s\subset \bar L$. Since the monodromy of $\cG_s$ is induced from $\cL_X$, it must be trivial. In other words, $\cG_s = k_{K_s}$. In conclusion, we get
$$0\rightarrow \cF_0 \rightarrow \cF \rightarrow \cL_X \rightarrow \bigoplus_{K_s\subset \bar L} i_{s*} k_{K_s} \rightarrow 0.$$
By Proposition \ref{strongredeq}, $\cF$ is not strongly reduced.

\smallskip

Suppose $\cF$ is not strongly reduced. By Proposition \ref{strongredeq} (2) and (3), there exists a (possibly empty) sublink $\bar{L}\subset L$ and an exact sequence $0\rightarrow \cF' \rightarrow \cF \rightarrow \cL_X \rightarrow \bigoplus_{K_s\subset \bar L} i_{s*} k_{K_s} \rightarrow 0$. Apply $j^{-1}$ to the exact sequence, we get
$$0\rightarrow j^{-1}\cF' \rightarrow j^{-1}\cF \rightarrow \cL_{X\setminus L} \rightarrow 0$$
This is also a short exact sequence of $\pi_L$ representations, where the first two terms as denoted by $(\rho', V')$, $(\rho,V)$. Because $\cL_{X\setminus L}$ has positive dimension, $(\rho',V')$ is a proper sub-representation of $(\rho,V)$. By Proposition \ref{stabsubuniv}, $(\rho_0,V_0)$ is a sub-representation of $(\rho',V')$. Hence, $(\rho_0,V_0)$ is proper sub-representation of $(\rho, V)$, and it is impossible to have $\cF = \cF_0$. Therefore, $\cF$ is not stable.

\end{proof}

\begin{rmk}\label{rsnsl}
For knots, being stable, reduced, or strongly reduced are equivalent. It allows us to reasonably classify simple sheaves \cite[Theorem 1.1]{Gao2}. The same classification for links will be harder.
\end{rmk}

\subsection{Sheaf moduli}\label{Sec:shmod}

Let $Sh^{s,0}_{\Lambda_L}(X)/Loc(X)$ be the subcategory of $Sh_{\Lambda_L}(X)/Loc(X)$ whose objects comes from $Sh^{s,0}_{\Lambda_L}(X)$.
Define
\begin{equation}\label{modshtil}
\widetilde{\cM}:= Ob(Sh^{s,0}_{\Lambda_L}(X)/Loc(X)).
\end{equation}
The equivalence relation in $\widetilde{\cM}$ is generated by $\cF_1\sim \cF_2$ if there exists a distinguished triangle
$$\cF_1\rightarrow \cF_2\rightarrow \cL \xrightarrow{+1},$$
where $\cF_1,\cF_2 \in Sh^{s,0}_{\Lambda_{L}}(X)$ and $\cL \in Loc(X)$.

\medskip

Let $\cF\in Sh(X)$. We denote $\cH^i:= H^i\cF$ to be the $i$-th cohomological sheaf. 

\begin{lem}\label{lemShclassify}

Any isomorphism class in $\widetilde{\cM}$ can be represented by an object $\cF$ with the following properties:
\begin{enumerate}
\item $\cH^i=0$ for $i\neq 0,1$; and
\item there is a decomposition $L = L_0\sqcup L_1$, such that $$\cH^1 \cong \bigoplus_{K_s\subset L_1} i_{s!}\cG_{\alpha_s},$$
where $\cG_{\alpha_s} \in loc(X)$ is a rank $1$ locally constant sheaf with monodromy $\alpha_s \in k^*$; and
\item $\cH^0$ is micro-supported along $\Lambda_{L_0}$.\end{enumerate}

\end{lem}

\begin{proof}

In this proof, we denote $\cF$ by $\cF^\bullet$ to emphasize that it is a chain complex of sheaves.

Let $x\in K_s\subset L$, and let $y \in B_x(\delta)\cap (X\setminus K_s)$. By setting $W^\bullet \cong (\cF^\bullet)_x, V^{\bullet}\cong (\cF^\bullet)_{y}$, the restriction map is a degree $0$ morphism $T^{\bullet}: W^\bullet \rightarrow V^{\bullet}$. Since $\cF^\bullet$ is simple, we have $cone(T^\bullet)\cong k$.

Because taking stalks is an exact functor, it intertwines with the kernel and cokernel functors, and therefore the cohomology functor. Hence, $\cH^i_x = (H^i\cF^\bullet)_x = H^iW^\bullet$, and $\cH^i_y = (H^i\cF^\bullet)_y = H^iV^\bullet$, and the local restriction map of $\cH^i$ is $T_{\cH^i} = H^iT^\bullet: H^iW^\bullet \rightarrow H^iV^\bullet$. Because the dg category of $k$-vector spaces is a semi-simple, each chain complex of vector spaces is equivalent to its cohomology complex with zero differentials. Passing to the mapping cone, there is a commutative diagram,
\begin{center}
\begin{tikzpicture}
  \node (A){$W^\bullet$};
  \node (B)[right of=A, node distance=3cm]{$V^\bullet$};
  \node (C)[below of=A, node distance=1.5cm]{$\oplus H^iW^\bullet[-i]$};
  \node (D)[right of=C, node distance=3cm]{$\oplus H^iV^\bullet[-i]$};
  \node (E)[right of=B, node distance=3.5cm]{$cone(T^\bullet)$};
  \node (F)[right of=E, node distance=2cm]{$$};
  \node (G)[right of=D, node distance=3.5cm]{$\oplus cone(T_{\cH^i}[-i])$};
  \node (H)[right of=G, node distance=2.5cm]{$.$};
  \draw[->] (A) to node [] {$$} (B);
  \draw[double] (B) to node [swap]{$$}(D);
  \draw[double] (A) to node [swap]{$$}(C);
  \draw[->] (C) to node []{$$}(D);
  \draw[->] (B) to node [] {$$} (E);
  \draw[->] (E) to node [] {$+1$} (F);
  \draw[->] (D) to node [] {$$} (G);
  \draw[->] (G) to node [] {$+1$} (H);
  \draw[->] (E) to node [] {} (G);
\end{tikzpicture}
\end{center}
The axiom of a triangulated category forces the third column to be a quasi-isomorphism. Note $cone(T_{\cH^i}) = \textrm{coker} (T_{\cH^i}) \oplus \ker (T_{\cH^i})[-1]$. We have
\begin{equation}\label{cohoTlocal}
H^i cone(T^\bullet) = \textrm{coker}(H^iW^\bullet \xrightarrow{T_{\cH^i}} H^iV^\bullet) \oplus \ker(H^{i+1}W^\bullet \xrightarrow{T_{\cH^{i+1}}} H^{i+1}V^\bullet).
\end{equation}

\noindent Since $cone(T^\bullet)\cong k$, there are two possibilities:
\begin{enumerate}
\item
$H^0W^\bullet \xrightarrow{T_{\cH^0}} H^0V^\bullet$ is injective with a rank $1$ cokernel, and $H^iW^\bullet \cong H^iV^\bullet$ for any $i\neq 0$, or
\item
$H^1W^\bullet \xrightarrow{T_{\cH^1}} H^1V^\bullet$ is surjective with a rank $1$ kernel, and $H^iW^\bullet \cong H^iV^\bullet$ for any $i\neq 1$.
\end{enumerate}
In either case, we have $H^iW^\bullet \cong H^iV^\bullet$ for $i \neq 0,1$. Therefore, $\cH^i \in loc(X)$ for $i \neq 0,1$.

Let $\tau_{<d}, \tau_{>d}, \tau_{\leq d}, \tau_{\geq d}$ be truncation functors. Apply $\tau_{<0}\rightarrow \textrm{id} \rightarrow \tau_{\geq 0}\xrightarrow{+1}$ to $\cF^{\bullet}$. Because $H^i(\tau_{<0}\cF) = H^i(\cF) =\cH^i$ if $i<0$ and $H^i(\tau_{<0}\cF) = 0$ if $i\geq 0$, we have $\tau_{<0}\cF \in Loc(X)$. Then we apply $\tau_{\leq 1} \rightarrow \textrm{id} \rightarrow \tau_{> 1} \xrightarrow{+1}$ to $\tau_{\geq 0}\cF^\bullet$. By a similar argument, $\tau_{>1} \cF^\bullet \in Loc(X)$. Therefore, 
$$\cF^\bullet \sim \tau_{\geq 0}\cF^\bullet \sim \tau_{\leq 1}\tau_{\geq 0}\cF^\bullet$$ in $\widetilde{\cM}$. This proves property (1).

We define a decomposition $L = L_0\sqcup L_1$, based on whether $cone(T_{\cH^1}) \cong 0$ or $cone(T_{\cH^0}) \cong 0$. In other words, $SS(\cH^0) \cap T^\infty X = \Lambda_{L_0}$, $SS(\cH^1) \cap T^\infty X = \Lambda_{L_1}$.

\smallskip

We consider $\cH^1$. For each $K_s\subset L_1$, there is a rank $1$ local system $\cG_{\alpha_s} \in loc(K_s)$ with monodromy $\alpha_s$, which is determined by the parallel transport of the microlocal Morse cone along the longitude. More specifically,  let $U_s$ be a small tubular neighborhood of $K_s$, let $\tilde{j}_s: U_s \rightarrow X$ be the open embedding, and let $i_{s_0}: K_s\rightarrow U_s$ be the closed embedding. Note $(i_{s_0})_* = (i_{s_0})_!$ and $\tilde{j}_s^{-1} = \tilde{j}_s^{!}$ Consider $\tilde{j}_s^{-1}\cH^1 = \tilde{j}_s^{-1}H^1\cF^\bullet = H^1(\tilde{j}_s^{-1}\cF^\bullet)\in Mod(U_s)$. Consider the short exact sequence $0\rightarrow k\rightarrow W_s\rightarrow V\rightarrow 0$ obtained from applying Lemma \ref{decompabsheaf}, we see a natural short exact sequence $0\rightarrow (i_{s_0})_!\cG_\alpha \rightarrow \tilde{j}_s^{!}\cH^1 \rightarrow \cL_{U_s}\rightarrow 0$. The morphism $(i_{s_0})_!\cG_\alpha \rightarrow \tilde{j}_s^{!}\cH^1$ can be extended globally, because 
$$
\textrm{Hom}_{U_s}((i_{s_0})_!\cG_{\alpha_s}, \tilde{j}_s^{!}\cH^1) 
	= \textrm{Hom}_{K_s}(\cG_{\alpha_s},i_{s_0}^! \tilde{j}_s^{!}\cH^1) 
	= \textrm{Hom}_{K_s}(\cG_{\alpha_s},i_s^{!}\cH^1) 	
	= \textrm{Hom}_{X}(i_{s!}\cG_{\alpha_s},\cH^1).
$$
Combining the morphisms obtained from all $K_s\subset L_1$, we get a morphism $\oplus_{K_s\subset L_1} i_{s!}\cG_{\alpha_s} \rightarrow \cH^1$. It can be completed to a short exact sequence $0\rightarrow \oplus_{K_s\subset L_1} i_{s!}\cG_{\alpha_s} \rightarrow \cH^1 \rightarrow \cL_X\rightarrow 0$. [The last term is a locally constant sheaf because of the micro-support condition.]

Consider the composition $\theta:  \cF^\bullet \rightarrow \cH^1[-1] \rightarrow \cL_X[-1]$, and complete it to a distinguished triangle
$$cone(\theta)[-1] \rightarrow \cF^\bullet \rightarrow \cL_X[-1]\xrightarrow{+1}.$$
Since $cone(\theta)[-1] \sim \cF^{\bullet}$ in $\cM$, we can assume $\cH^1 = \oplus_{K_s\subset L_1} i_{s!}\cG_{\alpha_s} $. This proves property (2).

Property (3) follows from the triangular inequality of micro-support.

\end{proof}

\begin{rmk}\label{lcscd}
The decomposition in Lemma \ref{lemShclassify} may not be unique. For example, $i_{s!}k_{K_S}[-1]$ and $j_{s!}k_{X\setminus K_s}$ live in different cohomological degrees, but they are isomorphic in $\widetilde{\cM}$ due to the following distinguished triangle
$$i_{s!}k_{K_S}[-1]\rightarrow j_{s!}k_{X\setminus K_s}\rightarrow k_X \xrightarrow{+1}. $$\hfill\qed
\end{rmk}

The sheaf $\cF$ in Lemma \ref{lemShclassify} can be described by extension classes in 
\begin{equation}\label{extdefM}
\mathrm{Ext}^1(\cH^1[-1], \cH^0).
\end{equation}

\noindent This extension class can be nontrivial. For each $i_{s!}\cG_{\alpha_s}$, we have
$$\mathrm{Ext}^1_X(i_{s!}\cG_{\alpha_s}[-1], \cH^0) = \mathrm{Ext}^2_{K_s}(\cG_{\alpha_s}, i^!\cH^0) = \mathrm{Ext}^2_{K_s}(\cG_{\alpha_s}, i^{-1}\cH^0[-2]) = \mathrm{Ext}^0_{K_s}(\cG_{\alpha_s}, i^{-1}\cH^0),$$
and the extension class is nonzero if $\cG_{\alpha_s^{-1}}\otimes i^{-1}\cH^0$ has nonzero global sections. Suppose the locally constant sheaf $i^{-1}\cH^0$ has rank $n$, its monodromy is given by an invertible matrix $M\in GL(n)$. If $\alpha_s$ is an eigenvalue of $M$, then the extension class in \eqref{extdefM} can be nontrivial.

\begin{defn}\label{shmoduli}
We define the subset $\cM\subset \widetilde{\cM}$, where an element $\cF\in \cM$ can be represented by $$\cF_{red}\oplus \cF_{deg}[-1],$$
where $\cF_{red} \in Sh^{s,0}_{\Lambda_{L_{red}}}(X)\cap Mod(X)$ is reduced, $\cF_{deg} = \oplus_{K_s\subset L_{deg}}i_{s!}\cG_{\alpha_s} \in Sh^{s,0}_{\Lambda_{L_{deg}}}(X)\cap Mod(X)$ for $\cG_{\alpha_s}\in loc(K_s)$ a rank $1$ local system with monodromy $\alpha_s$, and $L = L_{red} \sqcup L_{deg}$.
\end{defn}

In some cases, such as when $L$ is a knot or an unlink, we have $\cM = \widetilde{\cM}$. In general, we expect $\widetilde{\cM}$ is strictly larger than $\cM$. We give a criterion in Proposition \ref{cannotsimp} when the nontrivial extension class \eqref{extdefM} cannot be simplified.



\begin{lem} \label{lochomtrivial}
Let $K\subset X$ be a knot, and $i: K\rightarrow X$ the closed embedding. For $\alpha \in  k\setminus \{0,1\}$, let $\cG_\alpha \in loc(K)$ be the rank $1$ locally constant sheaf supported on $K$ with monodromy $\alpha$. For any $\cL \in Loc(X)$, 
$$\mathrm{RHom}(\cL,i_*\cG_\alpha) =0, \quad \mathrm{RHom}(i_!\cG_\alpha, \cL) =0.$$
\end{lem}
\begin{proof} Note $i_! = i_*$ and it is exact. Consider the special case $\cL = k_X$, then 
\begin{align}\label{vanishingwithloc}
\begin{split}
\mathrm{RHom}(\cL,i_*\cG_\alpha) &=\mathrm{RHom}(i^{-1}\cL, \cG_\alpha) = \mathrm{R}\Gamma(K, \cG_\alpha) =0. \\
\mathrm{RHom}(i_!\cG_\alpha, \cL) &= \mathrm{RHom}(\cG_\alpha, i^!\cL) = \mathrm{R}\Gamma(K, \cG_{\alpha^{-1}}[-2]) = 0.
\end{split}
\end{align}
Next we consider a general $\cL$. If $X = \bbR^3$, since the ambient space is homotopic to a point, then $\cL \cong \oplus_{n}k_X^{m_n}[d_n]$ is isomorphic to a direct sum of copies of $k_X$ with degree shifts. The vanishing follows from that in the special case. If $X = S^3$, choose $x\in X\setminus K$ and define $U:= X\setminus \{ x\}$. Note that $U$ is diffeomorphic to $\bbR^3$. Observe from \eqref{vanishingwithloc} that the desired $\mathrm{RHom}$ only depends on $\cL$ on a neighborhood of $K$. We take this neighborhood to be $U$, reducing the problem to the previous case when $X =\bbR^3$. This completes the proof.

\end{proof}

\begin{prop} \label{cannotsimp}
Suppose $\cF\in \widetilde{\cM}$ satisfy the properties in Lemma \ref{lemShclassify}. In addition, we assume:
\begin{enumerate}
\item $L_1$ is a knot. 
\item $\cG_\alpha \in loc(L_1)$ satisfies $\alpha \neq 1$.
\item $c\neq 0$ in the distinguished triangle
\begin{equation}\label{eq:nontriext}
\cH^0 \xrightarrow{a} \cF \xrightarrow{b} i_!\cG_\alpha[-1] \xrightarrow{c}\cH^0[1].
\end{equation}
\end{enumerate}
Then the distinguished triangle \eqref{eq:nontriext} does not split in $Sh(X)/Loc(X)$. In other words, there does not exist a morphism $p: i_!\cG_\alpha[-1] \rightarrow \cF$ in $Sh(X)/Loc(X)$, such that $b\circ p = \id$.
\end{prop}

\begin{proof}
We write $\cG := i_!\cG_\alpha[-1]$. A morphism $p: \cG\rightarrow \cF$ in the quotient category $Sh(X)/Loc(X)$ is represented by a roof $\cG \xleftarrow{u} \cG' \xrightarrow{v} \cF$, whereas $\cG' \in Sh(X)$, $u,v$ are morphisms in $Sh(X)$, and a distinguished triangle $\cG' \xrightarrow{u} \cG \rightarrow \cL\xrightarrow{+1}$ for some $\cL\in Loc(X)$. Since $\mathrm{RHom}(\cG,\cL)=0$ by Lemma \ref{lochomtrivial}, the distinguished triangle splits and we have $\cG' \cong \cG\oplus \cL[-1]$. We decompose $v: \cG' \rightarrow \cF$ as $v= (v_1,v_2)$. The assertion will follow from the claim that $b\circ v = 0$.

To prove the claim, we first observe that $b\circ v_2: \cL \rightarrow \cG$ vanishes due to Lemma \ref{lochomtrivial}. Then, we apply $\mathrm{RHom}(\cG, -)$ to \eqref{eq:nontriext} and the long exact sequence yields
$$\mathrm{Hom}(\cG, \cG[-1]) \rightarrow \mathrm{Hom}(\cG, \cH^0)\xrightarrow{g} \mathrm{Hom}(\cG, \cF)\rightarrow \mathrm{Hom}(\cG, \cG)\xrightarrow{h}\mathrm{Hom}(\cG, \cH^0[1])$$
Since $\mathrm{Hom}(\cG, \cG) = k$, and $h(1) = c \neq 0$ where $1 = \mathrm{id}_\cG$ and $c$ represents the nontrivial extension class, we have $h$ is an injection. Also note that $\mathrm{Hom}(\cG, \cG[-1]) = H^{-1}(S^1) =0$. Therefore $g$ is an isomorphism. In other words, any morphism $v_1: \cG\rightarrow \cF$ will factor through $\cH^0$, and can be expressed as $a \circ \tilde{v}_1$ for some $\tilde{v}_1: \cG\rightarrow \cH^0$. In the end, we find that $b\circ v_1 = b\circ a \circ \tilde{v}_1 = 0$.
\end{proof}

For an example satisfying the hypothesis of Proposition \ref{cannotsimp}, we consider the Hopf link $L = L_1 \sqcup L_2$, where each $L_i$ is an unknot. If $X = S^3$, then $X\setminus L_2$ homotopy retracts to $L_1$. Let $\cG$ be a rank $1$ local system on $L_1$ with monodromy $\alpha\neq 1$ and cohomological degree $1$, and $\cH$ be a rank $1$ local system on $X\setminus L_2$ with the same monodromy but cohomological degree $0$. Then the nontrivial class in $Ext^1(\cG,\cH)$ is nonsplit by Proposition \ref{cannotsimp}.

When nontrivial extension describe in Proposition \ref{cannotsimp} exists, the inclusion $\cM\subset \widetilde{\cM}$ is strict. For $L$ being a knot or an unlink, we have $\cM = \widetilde{\cM}$, but otherwise the nontrivial extension exists quite generally. Consider a link $L$ and a representation $(\rho, \pi_L)$ coming from a simple sheaf $\cF$. As long as there exists a loop $\gamma \in  \pi_L$ such that $\rho(\gamma)$ has a nontrivial eigenvalue, we can construct a nontrivial extension for the new link $L' := L \cup \gamma$, by gluing $\cF$ with a local system supported on $\gamma$.



\section{Correspondence}\label{Sec:correspondence}

We prove the main theorem in this section.

\begin{proof}[Proof of Theorem \ref{Mainthm}, part 1]
We first discuss the degenerate case. By Definition \ref{shmoduli}, an object $\cF\in \cM$ can be represented by $\cF_{red}\oplus \cF_{deg}[-1]$, where $\cF_{deg}[-1]$ is concentrated in cohomological degree $1$, and supported on a sublink $L_{deg}$. The correspondence augmentation $\epsilon$ has the property that, for any $K_s \in L_{deg}$, $\ep(\mu_s)=0$, $\ep(\lambda_s) = \alpha_s$ (defined in Definition \ref{shmoduli}), and $\ep(c_{ks}) =\ep(c_{sk}) =0$. In other words, all framed cords that start or end on $K_s$ are augmented to zero. With respect to this augmentation, the component $K_s$ is unlinked from the remaining components, which matches the sheaf description that the direct summand is a rank $1$ locally constant sheaf supported on $K_s$. The locally constant sheaf is parametrized by the monodromy, whereas the augmentation is parametrized by $\epsilon(\lambda_{s})$, and they match.
\end{proof}

From now on, we assume that $\cF$ is reduced. In particular, $T_s: W_s\rightarrow V$ (defined in Lemma \ref{decompabsheaf}) is injective for all $1\leq s\leq r$. We view $W_s$ as a subspace of $V$. Correspondingly, we assume that for any $K_s$, the augmentation $\ep$ does not vanish on all framed cords that start or end on $K_s$.

We shall construct the map between augmentations and sheaves. For each direction, the construction depends on some auxiliary data:
\begin{itemize}
\item From sheaves to augmentations, given a simple sheaf $\cF$, together with a set of local trivializations, we construct an augmentation $\ep_\cF$ in $\cA ug_{naive}$.
\item From augmentations to sheaves, given an augmentation, together with a braid representative of the link, we construct a sheaf $\cF_\ep$.
\end{itemize}
We will prove that these constructions are well-defined, and descend to a bijective correspondence between $\cA ug$ and $\cM$.

\subsection{From sheaves to augmentations}\label{Sec:fsta}


Suppose $\cF$ is a reduced sheaf. By Lemma \ref{decompabsheaf}, $\cF$ is equivalent to $(V, \rho, W_s, \rho_s, T_s)$.

\begin{defn}\label{defn:triv}
A \emph{local trivialization} for $\cF$, denoted by $f$, is an $r$-tuple $f = (f_1,\dotsb, f_r)$ of surjective linear transformations
$$f_s: V\rightarrow k,$$
such that $f_s|_{W_s} = 0.$
\end{defn}

For each $f_s$, we choose a right inverse $f_s^{-1}: k\rightarrow V$, and denote $f^{-1} :=(f_1^{-1},\dotsb, f_r^{-1})$.

\begin{defn}
Let $\cF\in Sh^{s,0}_{\Lambda_L}(X)\cap Mod(X)$ and $\cF \leftrightarrow (V,\rho,W_s,\rho_s,T_s)$. Let $f,f^{-1}$ be a local trivialization and its right inverse. Let $M_t: = \rho(m_t)$. Let $c_{st}$ be a framed cord from $K_s$ to $K_t$. Let $A_{c_{st}}$ be the trivialization map defined in Remark \ref{fungro}.

We define a map $\epsilon_{(\cF,f,f^{-1})}: \textrm{Cord}(L)\rightarrow k$ over its generators:
\begin{align}\label{ShtoAugdefn}
	\begin{split}
	(a) \quad &\epsilon_{(\cF,f,f^{-1})}(c_{st}) = f_s \circ A_{c_{st}} \circ (\idV - M_t)\circ {f}_t^{-1}, \\
	(b) \quad & \epsilon_{(\cF,f,f^{-1})}(\lambda_s) = {f}_s \circ A_{\ell_s} \circ {f}_s^{-1}, \\
	(c) \quad & \epsilon_{(\cF,f,f^{-1})}(\mu_s) = 1 - {f}_s \circ(\idV - M_s) \circ {f}_s^{-1}. \\ 
	\end{split}
\end{align}
\end{defn}

\begin{rmk}\label{calculateaug}
Suppose $n=\dim V$. We can choose a basis for $V$ and compute $\epsilon_\cF$ using matrix algebra. If vectors in $V$ are coorinatized as column $n$-tuples, then $f_s$ is a row vector, $f_s^{-1}$ is a column vector, $A_{c_{st}}$, $M_t$ are $n\times n$ matrix. Then \eqref{ShtoAugdefn} becomes matrix multiplications, ending with a $1\times 1$ matrix, i.e. an element in $k$.

Another way to compute $\epsilon_\cF$ is to take a vector  $v \in V\setminus W_t$. Then $f_t(v) \neq 0$ by Definition \ref{defn:triv}, and
\begin{equation}\label{calculateepsiloncij1}
\ep_{(\cF, f, f^{-1})}(c_{st}) = \frac{f_s \circ A_{c_{st}} \circ (\idV - M_t)(v)}{f_t(v)}.
\end{equation}
The result does not depend on the choice of $v$, because $v \cong v'$ in $V/W_t$ iff $v-v' = w\in W_t$, but $f_t(w) =0$, $(\idV - M_t)(w)=0$.

For a standard cord $\gamma_{ij}$, we have $A_{\gamma_{ij}} = \idV$. Let $v_j := (\idV - M_j)(v)$, then
\begin{equation}\label{indaugstdcord}
\ep_{(\cF, f, f^{-1})} (\gamma_{ij}) = \frac{f_i \circ (\idV - M_j)(v)}{f_j(v)} = \frac{f_i(v_j)}{f_j(v)}.
\end{equation}
\end{rmk}

\medskip
Next, we prove the following statements.
\begin{itemize}
\item
Fix $\cF, f $ and $f^{-1}$, then the map \eqref{ShtoAugdefn} defines an augmentation. 
\item
Fix $\cF$ and $f$, $\ep_{(\cF, f,f^{-1})}$ does not depend on the choice of $f^{-1}$. Therefore it makes sense to write $\ep_{(\cF, f)}$.
\item
Fix $\cF$, then two choice of trivializations define equivalent augmentations in $\cA ug$.
\end{itemize}
Together with Proposition \ref{prop:uptoloc}, these statements show that Definition \ref{ShtoAugdefn} yields a well-defined map from $\cM$ to $\cA ug$.

\begin{prop}
$\epsilon_{(\cF,f,f^{-1})}$ is an augmentation, i.e. it is a well-defined algebra morphism.
\end{prop}
\begin{proof}
For simplicity, we abbreviate $\epsilon_{(\cF,f,f^{-1})}$ by $\epsilon_{\cF}$ in this proof. We defined $\ep_\cF$ over generators, and we need to check that the relations in $\mathrm{Cord}(L)$ are preserved.

- Normalization. Let $[e_s]$ be the constant cord for $K_s$. We want to show
$$\epsilon_\cF(e_s) = 1 - \epsilon_\cF(\mu_s).$$
Because $A_{e_s} = \idV$, we have 
$\epsilon_\cF(e_s) =f_s \circ A_{e_s} \circ (\idV - M_s)\circ f_s^{-1} = f_s \circ (\idV - M_s)\circ f_s^{-1}.$ And by definition,
$1 - \epsilon_\cF(\mu_s) = f_s \circ(\idV - M_s) \circ f_s^{-1}.$

\smallskip
- Meridian. Suppose $c_{st}$ is a framed cord. We show $\epsilon_\cF(m_s\cdot c_{st}) = \epsilon_{\cF}(\mu_s)\epsilon_{\cF}(c_{st})$. The other $\epsilon_\cF(c_{st} \cdot m_t) = \epsilon_{\cF}(c_{st})\epsilon_{\cF}(\mu_t)$ is similar. The left hand side is
\begin{align*}
\epsilon_\cF(m_s\cdot c_{st}) &=  f_s \circ A_{m_s\cdot c_{st}} \circ (\idV - M_t)\circ f_t^{-1} \\
	&= f_s \circ M_s \circ A_{c_{st}} \circ (\idV - M_t)\circ f_t^{-1},
\end{align*}
and the right hand side is
\begin{align*}
& \epsilon_{\cF}(\mu_s)\epsilon_{\cF}(c_{st}) \\
	=&  \big(1 - f_s \circ(\idV - M_s) \circ f_s^{-1} \big) \big(f_s \circ A_{c_{st}} \circ (\idV - M_t)\circ f_t^{-1} \big)\\
	=& f_s \circ A_{c_{st}} \circ (\idV - M_t)\circ f_t^{-1} - f_s \circ(\idV - M_s) \circ f_s^{-1} \circ f_s  \circ A_{c_{st}} \circ (\idV - M_t)\circ f_t^{-1} \\
	=& f_s \circ A_{c_{st}} \circ (\idV - M_t)\circ f_t^{-1} - f_s \circ(\idV - M_s)  \circ A_{c_{st}} \circ (\idV - M_t)\circ f_t^{-1} \\
	=& f_s \circ M_s \circ A_{c_{st}} \circ (\idV - M_t)\circ f_t^{-1}.
\end{align*}
The third equality uses the following identity
\begin{equation}\label{compositionidentity}
(\mathrm{id}_V - M_s) \circ  f_s^{-1} \circ f_s = (\mathrm{id}_V - M_s).
\end{equation}
[For any $v\in V$, $f_s^{-1}\circ f_s(v) - v \in W_s$. Because $(\mathrm{id}_V - M_s)|_{W_s}=0$, the identity follows.]

\smallskip
- Longitude. We prove $\epsilon_\cF(\ell_s \cdot c_{st}) = \epsilon_{\cF}(\lambda_s)\epsilon_{\cF}(c_{st})$, and the argument for $\epsilon_\cF(c_{st} \cdot \ell_t) = \epsilon_{\cF}(c_{st})\epsilon_{\cF}(\lambda_t)$ is similar. The left hand side is
\begin{align*}
\epsilon_\cF(\ell_s \cdot c_{st}) &= f_s \circ A_{\ell_s \cdot c_{st}} \circ (\idV - M_t)\circ f_t^{-1} \\
	&=  f_s \circ A_{\ell_s} \circ A_{c_{st}} \circ (\idV - M_t)\circ f_t^{-1}.
\end{align*}
The right hand side is
\begin{align*}
\epsilon_{\cF}(\lambda_s)\epsilon_{\cF}(c_{st}) 
	&= (f_s \circ A_{\ell_s} \circ f_s^{-1})\circ (f_s  \circ A_{c_{st}} \circ (\idV - M_t)\circ f_t^{-1}).
\end{align*}
The two sides are equal, following:
$$f_s \circ A_{\ell_s} \circ f_s^{-1}\circ f_s = f_s \circ A_{\ell_s},$$
\noindent [because (1) $f_s^{-1}\circ f_s(v) -v \in W_s$, (2) $f_s \circ A_{\ell_s}\circ T_s = f_s \circ T_s \circ \rho_s(K_s)$, (3) $f_s \circ T_s =0$.]

\smallskip
- Skein relations. Let $c_{sk}$ and $c_{kt}$ be two composable cords. We need
\begin{equation}\label{shtoaugskein}
\epsilon_\cF(c_{sk} \cdot c_{kt}) - \epsilon_\cF (c_{sk}\cdot m_k\cdot c_{kt}) =   \epsilon_\cF(c_{sk})\epsilon_\cF(c_{kt}).
\end{equation}
The left hand side is
\begin{align*}
& \epsilon_\cF(c_{sk} \cdot c_{kt}) - \epsilon_\cF (c_{sk}\cdot m_k\cdot c_{kt}) \\
	=& \big( f_s \circ A_{c_{sk}\cdot c_{kt}} \circ (\idV - M_t)\circ f_t^{-1} \big) - \big( f_s \circ A_{c_{sk}\cdot m_k \cdot c_{kt}} \circ (\idV - M_t)\circ f_t^{-1} \big)\\
	=& \big( f_s \circ A_{c_{sk}} \circ A_{c_{kt}} \circ (\idV - M_t)\circ f_t^{-1} \big) - \big( f_s \circ A_{c_{sk}} \circ A_{m_k} \circ A_{c_{kt}} \circ (\idV - M_t)\circ f_t^{-1} \big)\\
	=&  f_s \circ A_{c_{sk}} \circ (\idV - A_{m_k})\circ A_{c_{kt}} \circ (\idV - M_t)\circ f_t^{-1} \\
	=& f_s \circ A_{c_{sk}} \circ (\idV - M_k)\circ A_{c_{kt}} \circ (\idV - M_t)\circ f_t^{-1}.
\end{align*}
And the right hand side is
\begin{align*}
\epsilon_\cF(c_{sk})\epsilon_\cF(c_{kt}) 
	=& \big(f_s \circ A_{c_{sk}} \circ (\idV - M_k)\circ f_k^{-1} \big) \big(f_k \circ A_{c_{kt}} \circ (\idV - M_t)\circ f_t^{-1} \big)\\
	=& f_s \circ A_{c_{sk}} \circ (\idV - M_k)\circ f_k^{-1} \circ f_k \circ A_{c_{kt}} \circ (\idV - M_t)\circ f_t^{-1} \\
	=& f_s \circ A_{c_{sk}} \circ (\idV - M_k)\circ  A_{c_{kt}} \circ (\idV - M_t)\circ f_t^{-1}.
\end{align*}
The second equality follows from \eqref{compositionidentity}.
\end{proof}

\begin{prop}
Fix $\cF$ and $f$, if $f^{-1}$ and $f'^{-1}$ are two right inverses, then $\ep_{(\cF,f,f^{-1})} = \ep_{(\cF,f,f'^{-1})}$.
\end{prop}
\begin{proof}
To check \eqref{ShtoAugdefn} (a) and (c), note that the image of $f_s^{-1} - f_s'^{-1}$ lies in the kernel of $f_s$, which is contained in $W_s$. By Lemma \ref{decompabsheaf} (3b), it vanishes when acted by $(\mathrm{id}_V - M_s)$. To check \eqref{ShtoAugdefn} (b), we apply Lemma \ref{decompabsheaf} (3a), and then use $f_s\circ T_s = 0$.
\end{proof}
\begin{prop}
Let $f = (f_1,\dotsb,f_r)$ and $ g= (g_1,\dotsb, g_r)$ be two local trivializations. Then $\epsilon_{(\cF,f)} \cong \epsilon_{(\cF,g)}$ in $\cA ug$.
\end{prop}
\begin{proof}
For any $v\in V\setminus W_s$, we define $d_s := {f_s(v)}/{g_s(v)}.$ Because $f_s(v)\neq 0$ and $g_s(v)\neq 0$, we have $d_s \in k^*$, and $(d_1,\dotsb, d_r)\in (k^*)^r$ is a dilation parameter. Note that $d_i$ does not depend on the choice of $v$. [If $v,v'\in V\setminus W_s$, let $a = f(v)/f(v')$, then $v- av' = \ker f = \ker g \Rightarrow a = g(v)/g(v') \Rightarrow f(v)/g(v) = f(v')/g(v')$.]

Suppose $c_{st}$ is a mixed cord from $K_s$ to $K_t$. Then
\begin{align*}
\epsilon_{(\cF,f)}(c_{st})
	&={f}_s \circ A_{c_{st}} \circ (\idV - M_t)\circ{f}_t^{-1} \\
	&= (d_s\cdot{g}_s) \circ A_{c_{st}} \circ (\idV - M_t)\circ ({g}_t^{-1}/{d_t}) 
	= (d_s/d_t) \cdot \epsilon_{(\cF,g)}(c_{st}).
\end{align*}
Therefore  $\epsilon_{(\cF,f)} \cong \epsilon_{(\cF,g)} $ in $\cA ug$.
\end{proof}

\begin{prop}\label{prop:uptoloc}
\smallskip
Suppose $\cF, \tilde{\cF} \in Sh^{s,0}_{\Lambda_{L}}(X)\cap Mod(X)$ are related by one of the following exact sequences,
\begin{enumerate}
\item
$0\rightarrow \tilde{\cF} \rightarrow \cF \rightarrow \cL_X \rightarrow 0$,
\item
$0\rightarrow \cL_X\rightarrow \cF \rightarrow \tilde{\cF} \rightarrow 0$,
\end{enumerate}
where $\cL_X \in loc(X)$, then a local trivialization $f= (f_1,\dotsb, f_r)$ of $\cF$ induce a local trivialization $\tilde{f} = (\tilde{f}_1,\dotsb, \tilde{f}_r)$ of $\tilde{\cF}$. Moreover, 
$$\epsilon_{(\cF, f)} = \epsilon_{(\tilde{\cF}, \tilde{f})}.$$
\end{prop}
\begin{proof} Suppose $\cF \leftrightarrow (\rho, V, \rho_s, W_s, T_s)$ and $\tilde{\cF} \leftrightarrow (\tilde{\rho}, \tilde{V}, \tilde{\rho}_s, \tilde{W}_s, \tilde{T}_s)$.

\smallskip
Case (1). $\tilde{\cF}$ is a subsheaf of $\cF$, hence $\tilde{V}\subset V$, $\tilde{W}_i = W_s \cap \tilde{V}$, and $\tilde{T_s} = T_s|_{W_s}$. Define 
$$\tilde{f}_s := f_s |_{\tilde{V}}: \tilde{V}\rightarrow k.$$
Since $\ker \tilde{f}_s = \tilde{V}\cap \ker f_s = \tilde{W}_s$, $\tilde{f} = (\tilde{f}_1,\dotsb, \tilde{f}_r)$ is a local trivialization of $\tilde{\cF}$. 

A local trivialization $f_s: V\rightarrow k$ induces an isomorphism $V/W_s\cong k$. We also have the natural isomorphism $\tilde{V} / \tilde{W}_s = \tilde{V}/ (\tilde{V}\cap W_s) \xrightarrow{} V/W_s$. Then the following commutative diagram gives $\epsilon_{{\cF}}(c_{st})= \epsilon_{\tilde{\cF}}(c_{st})$.

\begin{center}
\begin{tikzcd}[row sep=scriptsize, column sep=scriptsize]
& \tilde{W}_s \arrow[dl, hook] \arrow[rr, "\tilde{T}_s"]  & &  \tilde{V} \arrow[dl, hook] \arrow[dd, "\tilde{A}_{c_{st}}\circ (\textrm{id}_{\tilde{V}} -\tilde{M}_t)"  near start] \arrow[rr] & & \tilde{V}/\tilde{W}_s \arrow[dl, equal] \\ 

W_s \arrow[rr, crossing over, "T_s"] & & V \arrow[rr, crossing over] & & {V}/{W}_s \\

& \tilde{W}_t \arrow[dl, hook] \arrow[rr, "\tilde{T}_t" ', near start] & & \tilde{V} \arrow[dl, hook]  \arrow[rr] & & \tilde{V}/\tilde{W}_t \arrow[dl, equal] \arrow[from = uu, "{\epsilon_{(\tilde{\cF},\tilde{f})}(c_{st})}"]\\

W_t \arrow[rr, "T_t" '] & & V \arrow[from=uu, crossing over, "A_{c_{st}}\circ (\idV -M_t)" ', near start] \arrow[rr] & & V/W_t \arrow[from = uu, crossing over,  "{\epsilon_{({\cF},f)}(c_{st})}", near start]\\
\end{tikzcd}
\end{center}

\noindent The argument for $\mu_s$ and $\ell_s$ is similar.

\smallskip
Case (2). Since $\cL_X \subset \cF$, it can be trivialized to a vector space $V_0$, such that
$$V_0\subset \cap_{s=1}^r W_s\subset V.$$
The short exact sequence implies that
$$\tilde{V} = V/V_0,\quad \tilde{W}_s = W_s/V_0,$$
and that $\tilde{T}_s: \tilde{W}_s\rightarrow \tilde{V}$ is isomorphic to what $T_i$ induces between quotient spaces.

Since $V_0\subset W_i$, the linear function $f_s: V\rightarrow k$ vanishes on $V_0$, and it further induces a well-defined map $\tilde{f}_s: V/V_0\rightarrow k$. By construction, $\tilde{f}_s$ vanishes on $W_s/V_0\subset V/V_0$. Therefore $\tilde{f} = (\tilde{f}_1,\dotsb, \tilde{f}_r)$ is a local trivialization of $\tilde{\cF}$. The proof for $\epsilon_{(\cF, f)} = \epsilon_{(\tilde{\cF}, \tilde{f})}$ is similar.
\end{proof}

The formula for pure cords is simpler.

\begin{prop}\label{purecordformula}
Suppose $\cF\in Sh_{\Lambda_L}^{s,0}(X)\cap Mod (X)$. Suppose $c_s$ is a pure cord starting and ending on $K_s$. Then,
\begin{align}\label{simplifiedinducedaug}
	\begin{split}
	\epsilon_{\cF}(\lambda_s) &= \mathrm{tr}(\rho(\ell_s)) - \mathrm{tr}(\rho_s(K_s)),\\
	\epsilon_{\cF}(\mu_s) &= \mathrm{tr}(\mathrm{id}_V-\rho(m_s))+1,\\
	\epsilon_{\cF}(c_{s}) &= \mathrm{tr} (\rho(c_s) - \rho(m_s\cdot c_s)).
	\end{split}
\end{align}
In particular, the formula does not depend on the choice of the local trivialization.
\end{prop}

\begin{proof}

Following the construction in Remark \ref{calculateaug}, we can write $f_s$ as a $n\times 1$ row vector, $f_s^{-1}$ as a $1\times n$ column vector, and $A_{c_s}$ and $(\mathrm{id}_V - M_i)$ as $n\times n$-matrices. Since the composition
$$f_s \circ A_{c_s}\circ (\mathrm{id}_V - M_s) \circ f_s^{-1}$$
is a $1\times 1$ matrix, it also equals to its trace. Therefore we have
\begin{equation}\label{purecordsimplify1}
\epsilon_{\cF}(c_s) = \textrm{tr}(f_s \circ A_{c_s} \circ (\mathrm{id}_V - M_s) \circ f_s^{-1}) = \textrm{tr}(A_{c_s}\circ (\mathrm{id}_V - M_s) \circ f_s^{-1} \circ f_s) = \mathrm{tr}(A_{c_s}\circ (\mathrm{id}_V - M_s)).
\end{equation}

\noindent We assume $c_s$ is a loop after a cord homotopy. Then  \eqref{purecordsimplify1} becomes
$$\epsilon_\cF(c_s) = \textrm{tr} \big(\rho(c_s)(\idV - \rho(m_s)\big) = \textrm{tr} \big((\idV - \rho(m_s)\rho(c_s)\big) = \textrm{tr} \big(\rho(c_s) - \rho(m_s\cdot c_s)\big).$$

The argument for $\epsilon_\cF(\mu_s)$ is similar. For $\epsilon_\cF(\lambda_s)$, consider the following diagram:
\begin{center}
\begin{tikzpicture}
  \node (A){$W_s$};
  \node (B)[right of=A, node distance=2cm]{$V$};
  \node (C)[below of=A, node distance=1.7cm]{$W_s$};
  \node (D)[right of=C, node distance=2cm]{$V$};
  \node (E)[right of=B, node distance=2cm]{$k$};
  \node (F)[right of=E, node distance=2cm]{$0$};
  \node (G)[right of=D, node distance=2cm]{$k$};
  \node (H)[right of=G, node distance=2cm]{$0$};
  \node (X)[left of=A, node distance = 2cm]{$0$};
  \node (Y)[left of=C, node distance = 2cm]{$0$}; 
  \draw[->] (A) to node [] {$T_s$} (B);
  \draw[->] (B) to node [swap]{$\rho(\ell_s)$}(D);
  \draw[->] (C) to node [swap]{$T_s$}(D);
  \draw[->] (B) to node [] {$f_s$} (E);
  \draw[->] (E) to node [] {} (F);
  \draw[->] (D) to node [swap] {$f_s$} (G);
  \draw[->] (G) to node [] {} (H);
  \draw[->] (E) to node [swap] {$\epsilon_\cF(\lambda_s)$} (G);
  \draw[->] (A) to node [swap] {$\rho_s(K_s)$} (C);
  \draw[->] (X) to (A);
  \draw[->] (Y) to (C);
\end{tikzpicture}
\end{center}

\noindent Since $\ep_\cF(\lambda_s) = f_s\circ \rho(\ell_s)\circ f_s^{-1}$, the diagram on the right commutes. Hence, $\ep(\lambda_s) = \mathrm{tr}(\rho(\ell_s)) - \mathrm{tr}(\rho_s(K_s))$.
\end{proof}

\begin{rmk}\label{fpsak}
For knots, the induced augmentation is given by \eqref{simplifiedinducedaug}, matching that in \cite{Gao2}.
\end{rmk}

\begin{prop}\label{shtoaugproperty}
Let $\cF$ be a reduced sheaf and $\cF \leftrightarrow (V, \rho, W_s, \rho_s, T_s)$.
\begin{enumerate}
\item If $\rho(m_s) = \mathrm{id}_V$, then $\epsilon_\cF$ has the property that $R_j = 0$ for all $\{j\} = s$.

\item If $\cF$ is stable, then $\epsilon_\cF$ has the property that $R^i\neq 0$ for all $i$.
\end{enumerate}
\end{prop}

\begin{proof}
(1) Since $I - M_j = \mathrm{id}_V$, we have $\ep_\cF(\gamma_{ij}) = f_i \circ (I - M_j) \circ f_j^{-1} =0$.

(2) Let $v\in V\setminus \cup_{i=1}^n W_i$, and let $v_j := (\mathrm{id}_V- \rho (m_j))v$. By \eqref{indaugstdcord},
$$\ep_\cF(\gamma_{ij}) = \frac{f_i\circ (\mathrm{id}_V- \rho (m_j))(v)}{f_j(v)} = \frac{f_i(v_j)}{f_j(v)}.$$
We prove by contrapositive. Suppose $R^i = 0$, then $f_i(v_j)=0$ for all $j$, i.e. $v_j \in \ker f_i$ for all $j$. Hence $V_0 =\mathrm{Span}_k\{v_j\} \subset \ker f_i$. But $\ker f_i \subset V$ is a codimensional $1$ subspace, hence $\cF$ is not stable, a contradiction. 
\end{proof}



\subsection{From augmentations to sheaves}\label{Sec:fats}

Given $\ep$, the construction of the associated $\cF_\ep$ takes three steps: (1) augmentation representation, (2) augmentation subsheaf, (3) augmentation sheaf.

\subsubsection{Augmentation representation}
Let $L$ be an oriented link and $L'$ its Seifert framing. Let $\textrm{Cord}(L)$ be its framed cord algebra. Suppose $\epsilon: \textrm{Cord}(L)\rightarrow k$ is an augmentation. We choose an $n$-strand braid $B$ whose closure is $L$. In \cite{Gao3}, we defined the \textit{augmentation representation} $(\rho^{sub}_\epsilon, V^{sub}_\epsilon)$ \cite[Theorem-Definition 2.16]{Gao3} (which was denoted by $(\rho_\epsilon, V_\epsilon)$ in \cite{Gao3}), and it does not depend on the choice of the braid representative \cite[Theorem 2.18]{Gao3}.

Recall the construction. Let $D$ be the disk transverse to the braid. Let $\gamma_{ij}$, $1\leq i, j\leq n$ be the standard cords. Suppose $x_0 \in D$ is the base point for $X\setminus L$, and we fix a capping path $p_1$ from $x_1$ to $x_0$. Define $p_i := \gamma_{i1}\cdot p_1$. Let $h\in \pi_1(X\setminus L, x_0)$ be a based loop. Define $n\times n$ matrices:
$$
R = 
\begin{pmatrix}
\epsilon_{11} & \dotsb & \epsilon_{1n} \\
\vdots & \ddots &\vdots \\
\epsilon_{n1} & \dotsb & \epsilon_{nn}
\end{pmatrix}\qquad 
R^h = 
\begin{pmatrix}
\epsilon(p_1^{-1}\cdot h\cdot p_1) & \dotsb & \epsilon(p_1^{-1}\cdot h\cdot p_n) \\
\vdots & \ddots &\vdots \\
\epsilon(p_n^{-1}\cdot h\cdot p_1) & \dotsb & \epsilon(p_n^{-1}\cdot h\cdot p_n)
\end{pmatrix}
$$
Define $(\rho^{sub}_\epsilon, V^{sub}_\epsilon)$ to be
$$V_\epsilon^{sub}: = \textrm{Span}_k\{R_j\}_{1\leq j\leq n},$$
and 
$$\rho^{sub}_\epsilon(h)R_j := R^h_j.$$
In particular, the actions of meridian generators are given by
\begin{equation}\label{MeridianAction}
\rho_{\epsilon}(m_t)R_j = R_j -\epsilon(\gamma_{tj})R_{t}, \quad \rho_\epsilon(m_t^{-1})R_j = R_j + \mu_{\{t\}}^{-1}\epsilon(\gamma_{tj})R_t.
\end{equation}

\begin{prop}\label{diaugeqrep}
If $\ep_1 \cong \ep_2 $ in $ \cA ug$, then $(\rho_{\ep_1}^{sub},V_{\ep_1}^{sub}) \cong (\rho_{\ep_2}^{sub},V_{\ep_2}^{sub})$ as $\pi_L$-representations. 
\end{prop}

\begin{proof}
Let $(d_1,\dotsb, d_r) \in (k^*)^r$ be dilation parameters such that $\epsilon_1 \cong (d_1,\dotsb, d_r) \cdot \ep_2$. It means that $\epsilon_1(c_{ij}) = d_i/d_j \cdot \epsilon_2(c_{ij})$ for any mixed cord $[c_{ij}]$ from $x_i$ to $x_j$. Write $K$ as an $n$-strand braid closure. Let $R_1$ (resp. $R_2$) be the $n\times n$ matrix of augmented standard cords, i.e. $(R_{1})_{ij} = \ep_1(\gamma_{ij})$ (resp. $(R_{2})_{ij} = \ep_2(\gamma_{ij})$). Let $D = \textrm{diag}(d_{\{1\}},\dotsb, d_{\{n\}})$ be an $n\times n$ diagonal matrix. There is $R_1 = DR_2D^{-1}$. Moreover, for any based loop $h \in \pi_K$, we have $\tilde{R}_1^{h} = D\tilde{R}_2^{h}D^{-1}$. Hence $(\rho_{\ep_1},V_{\ep_1})$ and $(\rho_{\ep_2},V_{\ep_2})$ differ by a change of basis, and are isomorphic as $\pi_K$-representations.
\end{proof}

\subsubsection{Augmentation subsheaf} 
A row vector $R^i$ induces a linear map $f_i: V_\ep^{sub}\rightarrow k$ in the following way -- if $v= \sum_{j=1}^n a_jR_j \in V_\ep^{sub}$, then $f_i(v) = R^i\cdot v = \sum_{i=1}^n \ep_{ij}a_j$. [To see that $f_i$ is well-defined, let $V_\epsilon^{pre} = \oplus_{j=1}^n k[R_j]$ be the formal linear span, and note $V_\epsilon^{sub}$ is a quotient of $V_\epsilon^{pre}$. Each $R^i$ defines a linear map on $V_\epsilon^{pre}$ and it descents to a linear function on $V_\ep^{sub}$. For any $1\leq i\leq n$, and any linear relation $\sum_{j=1}^n a_i R_j = 0$, the $i$-th row is $\sum_{j=1}^n a_i \ep(\gamma_{ij})= 0$. In other words, $f_i(\sum_{j=1}^n a_i R_j)=0$.] If $R^i \neq 0$, then $f_i$ is surjective.

We will define an augmentation subsheaf $\cF_\ep^{sub}$, which micro-supported within $\Lambda_{\{I'\}}$. If $I'' =\emptyset$, then this set of linear maps, denoted by $f = (f_1,\dotsb, f_n)$, is a \textit{canonical local trivialization}.

\begin{defn}
Let $\ep: \textrm{Cord}(L)\rightarrow k$ be an augmentation, we define the \textit{augmentation subsheaf} $\cF_\ep^{sub}$ to be the sheaf associated to $(V_\ep^{sub}, \rho_\ep^{sub}, W_i, \rho_i, T_i)$, where
$$
W_i = \ker f_i, \quad T_i: W_i = \ker f_i \hookrightarrow V \;\textrm{ for }\; 1\leq i\leq n.
$$
\end{defn}

\begin{prop}\label{augsubsheafwelldef}
The augmentation subsheaf $\cF_\ep^{sub}$ is well-defined.
\end{prop}
\begin{proof}
By \cite[Theorem 1.1]{Gao3}, $(\rho_\ep, V_\ep)$ is a well-defined $\pi_L$-representation. It remains to check the properties of $W_i$ and $T_i$.

We first verify that $\rho_\ep^{sub}(m_i)|_{W_i} = \textrm{id}_{W_i}$, i.e. $W_i\subset \ker (\textrm{id}_{V_\ep^{sub}} - \rho_\ep^{sub}(m_i))$. Let $v= \sum_j a_j R_j$, then $v\in W_i = \ker f_i$ if and only if $\sum_j a_j \ep_{ij} = 0$. By \eqref{MeridianAction}, we have
$$\rho_\ep^{sub}(m_i) (\sum_j a_j R_j) = \sum_j a_j R_j - (\sum_j a_j\ep_{ij}) R_i = \sum_j a_j R_j.$$
Therefore $(\rho^{sub}_\ep(m_i) - \textrm{id}_{V})(v)= 0$ for any $v \in W_i$.

Next we check the compatibility. Let $g_i: = \ell_{i,\tau_B(i)}\cdot \gamma_{\tau_B(i),i}$. We claim that $f_{\tau_B(i)}= f_i \circ \rho_{\ep}^{sub}(g_i)$, if $\ell_{i,\tau_B(i)}$ does not contain $\ast_{\{i\}}$, and $f_{\tau_B(i)}= \ep(\lambda_{\{i\}}^{-1}) \cdot f_i \circ \rho_{\ep}^{sub}(g_i)$, if $\ell_{i,\tau_B(i)}$ contains $\ast_{\{i\}}$. It follows from the claim that $\rho_{\ep}^{sub}(g_i): \ker f_{\tau_B(i)} \stackrel{\sim}{\rightarrow} \ker f_i$. Since $\rho_{\ep}^{sub}(g_i) = A_{\ell_{i, \tau_B(i)}} A_{\gamma_{\tau_B(i),i}} = A_{\ell_{i, \tau_B(i)}}$, there is
$$A_{\ell_{i, \tau_B(i)}}: W_{\tau_B(i)} \stackrel{\sim}{\rightarrow} W_i.$$
We get the desired compatibility among $W_i$.

Proof of the claim:
If $\ell_{i,\tau_B(i)}$ does not contain $\ast_{\{i\}}$, for any $1\leq j \leq n$, there is
$$[\gamma_{\tau_B(i),j}] = [{\ell_{i, \tau_B(i)}} \cdot \gamma_{\tau_B(i),j}] = [g_i\cdot \gamma_{i,\tau_B(i)}\cdot \gamma_{\tau_B(i),j}] = [g_i\cdot \gamma_{ij}].$$
Apply $\ep$ to this identity, then the left hand side equals to $f_{\tau_B(i)}(R_j)$. We show that the right hand side equals to $f_i\circ \rho_\ep^{sub}(g_i) (R_j)$. By construction, $\rho^{sub}_\ep(g_i) (R_j) = \tilde{R}_j^{g_i}$. If $\tilde{R}_j^{g_i} = \sum_j a_jR_j$, then
$$f_i\circ \rho_\ep^{sub}(g_i) (R_j) = f_i (\tilde{R}_j^{g_i}) = f_i(\sum_j a_jR_j) = \sum_{j} a_j \epsilon_{ij}.$$
Note the $i$-th entry of the equation $\tilde{R}_j^{g_i} = \sum_j a_jR_j$ is $\epsilon(g_i\cdot \epsilon_{ij}) = \sum_j a_j \epsilon_{ij}$. Hence,
$$f_i\circ \rho^{sub}_\ep(g_i) (R_j) = \sum_j a_j \epsilon_{ij} = \ep(g_i\cdot \gamma_{ij}),$$
as desired.

If $\ell_{i,\tau_B(i)}$ contains $\ast_{\{i\}}$, the longitude relation becomes
$$
[\gamma_{\tau_B(i),j}] 
= \lambda_{\{i\}}^{-1}\cdot [{\ell_{i, \tau_B(i)}} \cdot \gamma_{\tau_B(i),j}] 
= \lambda_{\{i\}}^{-1}\cdot [g_i\cdot \gamma_{ij}],
$$
and the rest of the proof is similar.

\end{proof}

\begin{prop}\label{augsubXbraid}
$\cF_\ep^{sub}$ does not depend on the the braid representative.
\end{prop}
\begin{proof}
In \cite[Theorem 1.1]{Gao3}, we proved that $(V_\ep^{sub},\rho_\ep^{sub})$ is independent from the braid representative. It remains to check for $W_i$.

Note if $\rho_\ep^{sub}(m_i) \neq \id_{V_\ep^{sub}}$, then $W_i$ is uniquely determined by the action of the meridian, which is independent from the braid. We focus on $W_i$ with $\rho_\ep^{sub}(m_i) = \id_{V_\ep^{sub}}$. By construction, $W_i$ is determined by $f_i$.

For a conjugation such as $\tilde{B} = \sigma_s B\sigma_s^{-1}$, we can assume there is no marked point $\ast$ in the braiding region of $\sigma_i$ or $\sigma_i^{-1}$. It is straight forward to compute that $\tilde{f}^i  = f^i$ for $i\neq s, s+1$, $\tilde{f}^{s+1} = f^s$, and $\tilde{f}^s = f^{s+1}\circ \rho_\ep^{sub}(m_s^{-1})$. It is consist with the isomorphism of the augGao3mentation representation on the punctured disk, whereas $\tilde{m}_{i} = m_i, \tilde{m}_{s+1} = m_s$, and $\tilde{m}_s = m_s\cdot m_{s+1}\cdot m_s^{-1}$.

For stabilizations, we see immediately that $\tilde{f}_{n+1} = \tilde{f}_{n} = f_n$ after a negative stabilization (using matrix (2.19) \cite{Gao3}). For a positive stabilization, matrix (2.18) in \cite{Gao3} yields
$$\mu_{\{n\}}\tilde{f}_n = \tilde{f}_{n+1} .$$
Since we assumed $R_n = 0$, in particular there is $\ep(\gamma_{nn}) = 0$. Hence $1- \mu_{\{n\}} = 0$, or simply $\mu_{\{n\}} = 1$. Again we obtain $\tilde{f}_{n+1} = \tilde{f}_{n} = f_n$ for a positive stabilization.
\end{proof}

\begin{prop}\label{diaugeqpresh}
If $\ep_1\cong \ep_2$ in $\cA ug$, then $\cF_{\ep_1}^{sub} \cong \cF_{\ep_2}^{sub}$.
\end{prop}
\begin{proof}
If $\epsilon_1 \cong (d_1,\dotsb, d_r) \cdot \ep_2$, then $R_1 = DR_2D^{-1}$ for $D = \textrm{diag}(d_{\{1\}},\dotsb, d_{\{n\}})$. The diagonal matrix defines an isomorphism of the augmentation representation (Proposition \ref{diaugeqrep}) and subspaces $W_i$ transform accordingly. 
\end{proof}

\begin{prop}\label{augsubXglobalsec}
$\Gamma(\cF_\ep^{sub}) = 0$.
\end{prop}
\begin{proof}
In general, if $\cF \leftrightarrow (V, \rho, W_i, \rho_i, T_i)$, then $\Gamma(\cF) = \cap_{1\leq i\leq n} W_i$. 

For $\cF_\ep^{sub}$, $\Gamma(\cF_{\ep}^{sub}) = \cap_{1\leq i\leq n} W_i = \cap_{1\leq i\leq n} \ker f_i$. By construction, any vector in $\cap_{1\leq i\leq n} \ker f_i$ is a linear combination of $R_j$ which equals to zero, which is the zero vector in $V_\ep^{sub}$.
\end{proof}

\subsubsection{Augmentation sheaf}

We construct $\cF_\ep$. 

Recall the transverse disk $D$ and the index set $I = \{1, \dotsb, n\}$. Let $I = I' \cup I''$ such that $R^i \neq 0$ for $i\in I'$ and $R^i = 0$ for $i\in I''$. Define 
$$V_\ep := V_{\ep}^{sub }\oplus k.$$
Let $R_0$ be a basis vector for the direct summand $k$.

For $i\in I'$, define $W_i = W_i^{sub}\oplus k$. We have defined a surjective map $f_i: V_\ep^{sub} \rightarrow k$ via the row vector $R^i$. It can be naturally extended to a surjective linear map $f_i^{can}: V_\ep \rightarrow k$ such that $W_i: =\ker f_i^{can}$. The action of the meridian extends by identity on the direct summand $k$.

For $i\in I''$, let $M_i = \rho_\ep(m_i)$, define $M_i(R_j) = R_j$ for $1\leq j\leq n$ and $M_i(R_0) = R_0 + R_i$. Then $M_i$ is a unipotent linear transformation on $V_\ep$. such that $\mathrm{im}(\mathrm{id} - M_i)$ is spanned by $R_i$. We have $W_i = V_\ep^{sub}$ and $f_i^{can}$ is uniquely defined by $\ker f_i^{can} = W_i$ and $ (f_i^{can})^{-1} (1) = R_0$.

We obtain the augmentation sheaf, denoted by $\cF_\ep$, with canonical trivializations $f^{can}$.

The augmentation sheaf $\cF_\ep$ is an extension of the augmentation subsheaf $\cF_{\ep}^{sub}$. From the construction, we have an exact sequence of sheaves:
\begin{equation}\label{augsheafext}
0\rightarrow \cF_{\ep}^{sub}\rightarrow \cF_{\ep} \rightarrow k_X \rightarrow \bigoplus_{s\in I''}i_{s!}k_{K_s}\rightarrow 0,
\end{equation}
where $k_X $ is the constant sheaf on $X$. Let $\Lambda = \Lambda' \cup \Lambda''$ defined by the partition $\{1,\dotsb, r\} = \{I'\} \cup \{I''\}$. Then $SS(\cF_{\ep_\cF}^{sub}) = \Lambda'$ and $SS(\cF/\cF_{\ep_\cF}^{sub}) = \Lambda''$.

By Proposition \ref{augsubXglobalsec}, we have $\Gamma(\cF_{\ep}) = 0$.

\begin{prop}
The augmentation sheaf $\cF_\ep$ is well-defined.
\end{prop}
\begin{proof}
It suffices to check the compatibility of $W_i$ and $M_i: = \rho_\ep(m_i)$ for $i\in I''$. Since $W_i = V_\ep^{sub}$ for all $i\in I''$, they are compatible.

Define $h_j := \gamma_{j, \tau_B^{-1}(j)}\cdot \ell_{\tau_B^{-1}(j),j}$. We claim that: $\rho_\ep(h_i) R_j = A_{\ell_{\tau_B^{-1}(j),j}} R_j = R_{\tau_B^{-1}(j)}$, if $\ell_{\tau_B^{-1}(j),j}$ does not contain $\ast_{\{j\}}$, and $\rho_\ep(h_i) R_j = A_{\ell_{\tau_B^{-1}(j),j}} R_j =R_{\tau_B^{-1}(j)} \cdot \ep(\lambda_{\{i\}}^{-1})$, if $\ell_{\tau_B^{-1}(j),j}$ contains $\ast_{\{j\}}$. This gives the compatible actions of the meridians.

Proof of the claim: Note that $\rho_\ep|_{V_{\ep}^{sub}} = \rho_\ep^{sub}$. If $\ell_{\tau_B^{-1}(j),j}$ does not contain $\ast_{\{j\}}$. For any $1\leq i \leq n$, there is an identity of framed cords:
\begin{equation}\label{braidcordlemcol}
[\gamma_{i,\tau_B^{-1}(j)}] = [\gamma_{i,\tau_B^{-1}(j)} \cdot \ell_{\tau_B^{-1}(j), j}] = [\gamma_{i,\tau_B^{-1}(j)} \cdot \gamma_{\tau_B^{-1}(j), j} \cdot h_j] = [\gamma_{ij} \cdot h_j].
\end{equation}
After applying $\ep$ to the identity, the left hand side is the $i$-th entry of $R_{\tau_B^{-1}(j)}$, and the right hand side is the $i$-th entry of $\tilde{R}^{h_j}_j$, hence $\tilde{R}^{h_j}_j = R_{\tau_B^{-1}(j)}$. By the definition of the augmentation representation, there is $\rho_\ep(h_j)R_j = \tilde{R}^{h_j}_j$. Combing these equations, we have
$$\rho_\ep(h_i) R_j = R_{\tau_B^{-1}(j)}.$$

If $\ell_{\tau_B^{-1}(j),j}$ contains $\ast_{\{j\}}$, then $[\gamma_{i,\tau_B^{-1}(j)}] = [\gamma_{ij} \cdot h_j]\cdot \lambda_{\{j\}}^{-1}$, and the rest is similar.
\end{proof}

\begin{prop}
$\cF_\ep$ does not depend on the braid reprentative.
\end{prop} 
\begin{proof}
By Proposition \ref{augsubXbraid}, $\cF_{\ep}^{sub}$ does not depend on the braid. The extension class \eqref{augsheafext} does not depend on the braid either. 
\end{proof}

\begin{prop}
If $\ep_1 \cong \ep_2$ in $\cA ug$, then $\cF_{\ep_1} \cong \cF_{\ep_2}$.
\end{prop}
\begin{proof}
Both $\cF_{\ep_1}$ and $\cF_{\ep_2}$ are the same extension class in \eqref{augsheafext}, hence they are equal.
\end{proof}

\begin{prop}
Let $(\rho_\ep^{sub}, V_\ep^{sub}) \leftrightarrow \cE_\ep^{sub} \in loc(X\setminus K)$, and let $j: X\setminus K \rightarrow X$ be the open embedding. If $\ep$ is generic, then
$$\cF_\ep = \cF_\ep^{sub} = j_*\cE_\ep.$$
\end{prop}
\begin{proof}
Because $I'' =0$, there is no extension, hence $\cF_\ep = \cF_\ep^{sub}$. Because $J'' = 0$, $\rho_\ep(m_i) \neq \id_{V_\ep}$, and $W_i$ is uniquely determined by the invariant subspace of $\rho_\ep(m_i)$, i.e. $\ker (\textrm{id}_{V_\ep} - \rho_{\ep}(m_i))$.
\end{proof}

\begin{prop}\label{augtoshproperty}
Let $\ep$ be an augmentation. 
\begin{enumerate}
\item If $R_j = 0$, then $\rho_\ep(m_j) =\mathrm{id}_{V_\ep}$.
\item If $R^i\neq 0$ for all $i$, then $\cF_\ep$ is stable.
\end{enumerate}
\end{prop}
\begin{proof}
(1) By \eqref{MeridianAction}, $\rho_\ep^{sub}(m_j) = \mathrm{id}_{V_\ep^{sub}}$, and hence $\rho_\ep(m_j) =\mathrm{id}_{V_\ep}$ by construction.

(2) Under the hypothesis, we have $\cF_\ep = \cF_\ep^{sub}$. Then by construction $V_\ep$ is spanned by the image of $(\id_{V_\ep} - \rho_\ep(m_j))$ for all $j$.
\end{proof}

\subsection{Proof of the correspondence.}

Following Proposition \ref{shtoaugproperty} and \ref{augtoshproperty}, we see that generic augmentations correspond to sheaves that are stable and $\rho(m)\neq \mathrm{id}$ for any meridian $m$. We say such sheaves are \textit{generic}. By the same propositions, we have the following corresponding properties between sheaves and augmentations. (Recall the index sets from Definition \ref{defindexsets}.)

\begin{table}[!htbp]
	\setlength{\tabcolsep}{20pt}
	\renewcommand*{\arraystretch}{1.3}
\begin{tabular}{|c|c|}
\hline
sheaves          & augmentations \\ \hline
reduced          &       $I''\cap J''=\emptyset$       \\ \hline
stable           &       $I''=\emptyset$       \\ \hline
$\rho(m)\neq \id$ &   $J''=\emptyset$           \\ \hline
generic          &     $I'' = J'' = \emptyset$   \\ \hline
\end{tabular}
\end{table}

\begin{proof}[Proof of Theorem \ref{Mainthm}, part 2] The theorem follows from the two statements:

\begin{enumerate}
\item If $I''\cap J'' = \emptyset$ for $\ep$, then $\ep_{(\cF_\ep, f^{can})} = \ep$.
\item Suppose $\cF$ is reduced, then $\cF_{\ep_\cF} \sim \cF$.
\end{enumerate}

Proof of (1). For simplicity, denote $\ep'= \ep_{(\cF_\ep, f^{can})}$ and $f_i = f_i^{can}$. We first show $\ep_{ij}' = \ep_{ij}$ for generic entries. Take $v\in V\setminus W_j$, by \eqref{calculateepsiloncij1},
$$\ep_{ij}' = \frac{f_i \circ A_{\gamma_{ij}} \circ (\idV - M_j)(v)}{f_j(v)} = \frac{f_i \circ (\idV - M_j)(v)}{f_j(v)}.$$
Suppose $v = \sum a_sR_s$. By Lemma \ref{MeridianAction}, we have $(\idV - M_j)(R_s) = \ep_{js}R_j$. Also recall the canonical $f_i$ is defined from the row vectors of the matrix $R$. Hence, 
$$v \stackrel{\idV - M_j}{\longmapsto} \sum a_s (\ep_{js} R_j) = (\sum a_s\ep_{js}) R_j \stackrel{f_i}{\longmapsto} \ep_{ij} (\sum a_s\ep_{js}),$$
and $f_j(v) = \sum a_s\ep_{js} \neq 0$, [because $v\notin W_j$]. Taking the ratio, we get $\ep'_{ij} = \ep_{ij}$.

For non-generic entries, note $R^i$ and $R_i$ cannot be simultaneously zero, $R^j$ and $R_j$ cannot be simultaneously zero

If $R^i =0$, then $V_\ep^{sub}\subset \ker f_i$, and hence $\ep_{ij}' = f_i(v_j)/ f_j(v) =0 = \ep_{ij}$. 

If $R_j=0$, then $M_j = \mathrm{id}$, and $\ep_{ij}' = 0 = \ep_{ij}$. 

If $R_i= 0$ and $R^j = 0$, note it implies that $R^i \neq 0$. Recall we have chosen $f_j^{-1}$ such that $(\mathrm{id} - M_j) \circ f_j^{-1}(1) = R_j$. Then $\ep_{ij}' = f_i(R_j) = \ep_{ij}$.

\medskip

Proof of (2). We first prove that $\cF_0 \cong \cF_{\ep_{\cF}}^{sub}.$

Take an $n$-strand braid $B$. Let $\cF \leftrightarrow (V, \rho, W_i, \rho_i, T_i)$, $1\leq i\leq n$. Take $v\in V\setminus \cup_{i=1}^n W_i$ and choose a trivialization $f$. Define $v_i = (\mathrm{id}_V - M_i)v$. Note $v_i\neq 0$. By Remark \ref{calculateaug}, we have $$\ep_\cF(\gamma_{ij}) = f_i(v_j)/f_j(v).$$

We define a morphism $V_{\ep_\cF}\rightarrow V$ by 
\begin{equation}\label{pcp2ce}
R_i\mapsto v_i/f_i(v).
\end{equation}
By definition, the image is $V_0$. To see it is well-defined, suppose $\sum_j a_j R_j =0$, then for any $i$, $\sum_j a_j\ep_\cF (\gamma_{ij}) = 0 \Leftrightarrow \sum_j a_j f_i(v_j)/f_j(v)= 0 \Leftrightarrow f_i(\sum a_j v_j/f_j(v)) =0$. Next, since $\Gamma(\cF)=0$, we have $\cap_{i=1}^{n} \ker f_i = \{0\}.$ Therefore $\sum_j a_j R_j =0$ if and only if $\sum_j a_j v_j/f_j(v) =0$. We conclude $V_{\ep_\cF}\cong V_0$. Next we check the actions of meridians. Since $\ker (\mathrm{id}_V - M_i) = \ker f_i$, we have $(\mathrm{id}_V - M_i) v_j = f_i(v_j)/f_i(v) \cdot (\mathrm{id}_V - M_i)v = f_i(v_j)/f_i(v)\cdot v_i$. By \eqref{MeridianAction}, we have
$(\mathrm{id}-\rho_{\epsilon_\cF}(m_i))R_j = \epsilon_\cF(\gamma_{ij})R_{i}$. Invoking \eqref{pcp2ce}, we have
$$(\mathrm{id}_V - M_i) v_j/f_j(v) = f_i(v_j)/f_i(v)\cdot v_i/f_j(v) = \epsilon_\cF(\gamma_{ij}) \cdot v_i/f_i(v).$$
The subspaces $W_i$ and the maps $T_i$ are naturally induced. Therefore, $\cF_0 \cong \cF_{\ep_{\cF}}^{sub}.$

Next we prove that $\cF_{\ep_\cF} \sim \cF$. Let $\{I''\}$ be the image of $\{-\}|_{I''}: I''\rightarrow \{1,\dotsb, r\}$. We have an exact sequence
\begin{equation}\label{corrext1}
0\rightarrow \cF_0\rightarrow \cF \rightarrow \cL_X \rightarrow \bigoplus_{s\in I''} i_{s!}k_{K_s}\rightarrow 0,
\end{equation}
for some (locally) constant sheaf $\cL_X = k_X^n\in loc (X)$ of rank $n$. Also, applying \eqref{augsheafext} to $\cF_{\ep_{\cF}}$, we obtain
\begin{equation} \label{corrext2}
0\rightarrow \cF_{\ep_\cF}^{sub}\rightarrow \cF_{\ep_{\cF}} \rightarrow k_X\rightarrow \bigoplus_{s\in I''}i_{s!}k_{K_s} \rightarrow 0,
\end{equation}
where $k_X $ is the constant sheaf on $X$. Since $\cF_0 \cong \cF_{\ep_\cF}^{sub}$, both exact sequences \eqref{corrext1} and \eqref{corrext2} arise from $\mathrm{Ext}^2(\oplus_s i_{s!}k_{K_s},\cF_0) = \oplus_s \mathrm{Ext}^2(i_{s!}k_{K_s},\cF_0) = \oplus_s \mathrm{Ext}^2(k_{K_s},i_{s}^{-1}\cF_0)$, where the direct sum is over $I''$. For a fixed component $K_s$, $\mathrm{Ext}^2(k_{K_s},i_{s}^{!}\cF_0) = R^0\Gamma(K_s, i_{s}^{!}\cF_0[2])$, which only depends on $\cF_0$ restricted to a tubular neighborhood $U_s$ of $K_s$. Let $m_s$ be a meridian of $K_s$, then $\rho(m_s)$ is a uniponent matrix (following the construction of $\cF_0$), and so is $\rho_{\ep_\cF}(m_s)$ (following the construction of $\ep \mapsto \cF_\ep$). Now, we restrict both \eqref{corrext1} and \eqref{corrext2} to $U_s$, and we get
\begin{equation}\label{corrext1a}
0\rightarrow \cF_0|_{U_s} \rightarrow \cF|_{U_s} \rightarrow k_{U_s\setminus K_s} \oplus (k_{U_s})^{n-1} \rightarrow 0,
\end{equation}
and 
\begin{equation}\label{corrext2a}
0\rightarrow \cF_{\ep_\cF}^{sub}|_{U_s} \rightarrow \cF_{\ep_{\cF}}|_{U_s} \rightarrow k_{U_s\setminus K_s} \rightarrow 0.
\end{equation}
Here $n = \textrm{rank } \cL_X$, and we have replaced $\cL_X|_{U_s} \rightarrow i_{s!}k_{K_s}$ and $k_{U_s} \rightarrow i_{s!}k_{K_s}$ by their kernels. Comparing $\rho(m_s)$ and $\rho_{\ep_\cF}(m_s)$ (namely the kernels and images of $\id_{V_0} - \rho(m_s)$ and $\id_{V_{\ep_{\cF}}} - \rho_{\ep_\cF}(m_s)$), we find that \eqref{corrext1a} splits into a direct sum of \eqref{corrext2a} and $0 \rightarrow 0 \rightarrow (k_{U_s})^{n-1}\rightarrow (k_{U_s})^{n-1}\rightarrow 0$. Therefore they are given by the same extension class in $\mathrm{Ext}^2( i_{s!}k_{K_s},\cF_0)$. Collecting all indices $s\in I''$, we deduce that \eqref{corrext1} and \eqref{corrext2} arises from the same class in $\mathrm{Ext}^2(\oplus_s i_{s!}k_{K_s},\cF_0)$. It yields that,
$$cone(\cF \rightarrow \cL_X) \cong cone(\cF_{\ep_{\cF}} \rightarrow k_X).$$
Therefore, we conclude that $\cF \sim \cF_{\ep_\cF}$ in $\cM$ (i.e. they are equivalent up local local systems).
\end{proof}

\subsection{An example} Let $L$ be the $3$-component unlink, and let $B\in Br_3$ be the empty word. Consider the following augmentation:
$$
R=
\begin{pmatrix}
0 & \epsilon_{12} & \epsilon_{13} \\
0 & 0 & 0 \\
0 & \epsilon_{32} & \epsilon_{33}
\end{pmatrix}.
$$
We assume that $\ep_{12}$, $\ep_{13}$, $\ep_{32}$, $\ep_{33}$ and 
$
\det \begin{pmatrix}
\ep_{12} & \ep_{13} \\ \ep_{32} & \ep_{33}
\end{pmatrix}
$ are nonzero.

In the augmentation representation $(V_{\ep}^{sub}, \rho_{\ep}^{sub})$, we have $V_\ep = \textrm{Span}_k\{R_2,R_3\} \cong k^2$, and with respect to the ordered basis,
$$
\rho_\ep^{sub}(m_1) = \begin{pmatrix}
1 & 0 \\ 0 &1
\end{pmatrix},
\quad
\rho_\ep^{sub}(m_2) = \begin{pmatrix}
1 & 0 \\ 0 &1
\end{pmatrix},
\quad
\rho_\ep^{sub}(m_3) = \begin{pmatrix}
1 & 0 \\ -\ep_{32} &1- \ep_{33}
\end{pmatrix}.
$$
The augmentation sheaf $\cF_\ep$ is given by
$$V_\ep = k^3, 
\quad 
M_1= \begin{pmatrix}
1&& \\
&1& \\
&&1
\end{pmatrix},
\quad 
M_2= \begin{pmatrix}
1&& \\
1&1& \\
&&1
\end{pmatrix},
\quad 
M_3=\begin{pmatrix}
1&& \\
&1& \\
&-\ep_{32}&1-\ep_{33}
\end{pmatrix},
$$
where $M_j = \rho_\ep(m_j)$, and
$$
W_1 = \textrm{Span}_k 
\bigg\{\begin{pmatrix}
1 \\ 0 \\ 0
\end{pmatrix}, 
\begin{pmatrix}
0 \\ -\ep_{13} \\ \ep_{12}
\end{pmatrix}
\bigg\}
,
\quad
W_2 = \textrm{Span}_k 
\bigg\{\begin{pmatrix}
0 \\ 1 \\ 0
\end{pmatrix}, 
\begin{pmatrix}
0 \\ 0 \\ 1
\end{pmatrix}
\bigg\}
,
\quad
W_3 = \textrm{Span}_k 
\bigg\{\begin{pmatrix}
1 \\ 0 \\ 0
\end{pmatrix}, 
\begin{pmatrix}
0 \\ -\ep_{33} \\ \ep_{32}
\end{pmatrix}
\bigg\}
$$
The canonical trivializations ($f_i := f_i^{can}$) are
$$
f_1 = (0, \ep_{12},\ep_{13}),
\quad
f_2 = (1, 0, 0),
\quad
f_3 = (0, \ep_{32}, \ep_{33}),
$$
and we choose their right inverses to be
$$
f_1^{-1} = \begin{pmatrix}
0 \\ 0 \\ \ep_{13}^{-1}
\end{pmatrix},
\quad
f_2^{-1} = \begin{pmatrix}
1 \\ 0 \\ 0
\end{pmatrix},
\quad
f_3^{-1} = \begin{pmatrix}
0 \\ 0 \\ \ep_{33}^{-1}
\end{pmatrix}.
$$
Using the formula, $\ep_\cF(\gamma_{ij}) = f_i \circ (\id - M_j) \circ f_j^{-1}$, it is straightforward to verify that $\ep_\cF(\gamma_{ij}) = \ep_{ij}$.

\end{document}